\documentclass[a4paper]{article}
\usepackage{amsmath,amssymb,amsfonts}
\usepackage{colordvi}
\usepackage[all]{xy}

\def\sur{\overline}
\def\sou{\underline}

\def\mpn{\medskip\par\noindent}

\def\mmpn{\vskip 1em minus 1em\par\noindent}

\def\smp{\smallskip\par}

\newcommand{\flh}[2]{\mathop{\hbox to  4ex{\rightarrowfill}}_{#2}^{#1}
\limits}

\def\Id{{\rm id}}
\def\Res{{\rm Res}}
\def\Ind{{\rm Ind}}
\def\Hom{{\rm Hom}}
\def\End{{\rm End}}

\def\Out{{\rm Out}}
\def\Mod{{\rm Mod}}
\def\Inf{{\rm Inf}}
\def\Def{{\rm Def}}

\def\Iso{{\rm Iso}}
\def\Conj{{\rm Conj}}

\def\Ker{{\rm Ker}}

\def\Defres{{\rm Defres}}

\def\Indinf{{\rm Indinf}}

\def\rank{{\rm rank}}

\def\op{^{op}}

\def\ls#1#2{{\,^{#1}\!#2}}

\def\Z{\mathbb{Z}}

\def\Q{\mathbb{Q}}

\def\pf{\par\bigskip\noindent{\bf Proof~: }}
\def\endpf{\nolinebreak~\leaders\hbox to 1em{\hss\
\hss}\hfill~\raisebox{.5ex}
{\framebox[1ex]{}}\par\bigskip}

\makeatletter
\renewenvironment{enumerate}{\ifnum \@enumdepth >3 \@toodeep\else
       \advance\@enumdepth \@ne
       \edef\@enumctr{enum\romannumeral\the\@enumdepth}\list
       {\csname  label\@enumctr\endcsname}{\setlength{\topsep}{1ex}
\setlength{\itemsep}{0 pt}\usecounter
         {\@enumctr}\def\makelabel##1{\hss\llap{##1}}}\fi}{\endlist}

\def\@seccntformat#1{\csname the#1\endcsname.\quad}

\def\section{\pagebreak[3]\setcounter{prop}{0}\setcounter{equation}{0}\@startsection{section}{1}{\z@}{4ex plus  6ex}{2ex}{\center\reset@font \large\bf}}

\makeatother

\def\theprop{\thesection.\arabic{prop}}

\newenvironment{enonce}[1]{\pagebreak[3]\refstepcounter{prop}\mmpn
{{\bf  \thesection.\arabic{prop}.\ #1.}}\begin{it} }{\end{it}\smp}
\def\thesection{\arabic{section}}

\newcommand{\result}[1]{\begin{enonce}{#1}}
\newcommand{\fresult}{\end{enonce}}

\begin{document}

\centerline{\Large\bf Simple biset functors and double Burnside ring}
\vspace{.5cm}
\centerline{Serge Bouc, Radu Stancu, and Jacques Th\'evenaz}
\vspace{1cm}

\begin{footnotesize}
{\bf Abstract~:} Let $G$ be a finite group and let $k$ be a field.
Our purpose is to investigate the simple modules for the double Burnside ring $kB(G,G)$. It turns out that they are evaluations at $G$ of simple biset functors.
For a fixed finite group $H$, we introduce a suitable bilinear form on $kB(G,H)$ and we prove that the quotient of $kB(-,H)$ by the radical of the bilinear form is a semi-simple functor. This allows for a description of the evaluation of simple functors, hence of simple modules for the double Burnside ring.

\bigskip\par

{\bf AMS Subject Classification~:} 19A22, 20C20.\par
{\bf Key words~:} biset, Burnside ring, simple module. 
\end{footnotesize}


\section{Introduction}
\noindent
Let $G$ be a finite group. The double Burnside ring $B(G,G)$ of all $(G,G)$-bisets plays a crucial role in many recent developments of representation theory and homotopy theory. It  appears in the theory of biset functors \cite{Bo3}, which was successfully used for the solution of important problems in representation theory (e.g. the classification of endo-permutation modules for a $p$-group \cite{Bo2}). In homotopy theory, a subring of $B(G,G)$ appears in the theory of $p$-completed classifying spaces \cite{MP,BF}, fusion systems and $p$-local groups \cite{BLO}, idempotents associated to fusion systems \cite{Ra, RS}.\par

One of the main issues concerning the ring structure is to understand the simple $B(G,G)$-modules,
since they appear in the semi-simple quotient of~$B(G,G)$ by its Jacobson radical. For this, it suffices to work over a field~$k$, and we consider the finite-dimensional $k$-algebra $kB(G,G)=k\otimes_{\Z}B(G,G)$.
Very little is known in general about the ring structure of $kB(G,G)$. It is known that $kB(G,G)$ is semi-simple only for cyclic groups in suitable characteristic, e.g.~in characteristic zero (see Section~6.1 in~\cite{Bo3}). For a $p$-group $P$ and for the subring of $kB(P,P)$ spanned by bisets which are free on one side, Benson and Feshbach~\cite{BF} have a description of all simple modules, while Henn and Priddy give a sufficient condition for this subring to be a local ring (and they prove that this condition occurs quite frequently).
More recently, for any finite group $G$, Boltje and Danz~\cite{BD} define a ghost ring for the subring of $kB(G,G)$ spanned by bisets which are free on one side (respectively on both sides) and this sheds some new light about the ring structure, especially over a field $k$ of characteristic zero. \par

The main purpose of this paper is to study the full double Burnside ring $kB(G,G)$ and analyze the simple $kB(G,G)$-modules, using their connection with simple biset functors. This connection is quite deep since any simple $kB(G,G)$-module determines uniquely a simple biset functor, and conversely any ($k$-linear) simple biset functor has an evaluation at $G$ which is a simple $kB(G,G)$-module (provided it is non-zero). The subring of $kB(G,G)$ generated by the bisets which are free on both sides is also considered briefly.\par

We work mainly over an algebraically closed field $k$ and consider the finite-dimensional $k$-vector space $kB(G,H)=k\otimes_{\Z}B(G,H)$, where $G$ and $H$ are finite groups and $B(G,H)$ is the Grothendieck group of $(G,H)$-bisets. Let
$$k\sur B(G,H)=kB(G,H)/kI(G,H) \,,$$
where $kI(G,H)$ is the $k$-subspace generated by all bisets which factor through a proper subquotient of~$H$. We define a canonical bilinear form on $k\sur B(G,H)$ and pass to the quotient by the radical $R(G,H)$ of this form.
Allowing $G$ to vary, we obtain a biset functor $k\sur B(-,H)$ and a subfunctor $R(-,H)$. We prove that $k\sur B(-,H)/R(-,H)$ is semi-simple, more precisely the largest semi-simple quotient of the biset functor $k\sur B(-,H)$.\par

Evaluation at~$G$ gives rise to a semi-simple $kB(G,G)$-module
$$k\sur B(G,H)/R(G,H) \,,$$
with an explicit decomposition into simple summands. This provides the main tool for analyzing simple $kB(G,G)$-modules, or equivalently, evaluation of simple biset functors. In particular, we obtain a formula for the dimension of the evaluation $S_{H,V}(G)$ of a simple biset functor $S_{H,V}$, generalizing the formula obtained in~\cite{Bo1} for the case where $V$ the trivial module.\par

After an introductory Section~2, we review in Section~3 the connections between simple modules for $kB(G,G)$ and simple functors. In Section~4, we define the standard quotient $k\sur B(G,H)$ of $kB(G,H)$ and prove some of its properties. We prove in Section~5 that any simple $kB(G,G)$-module has a minimal group attached to it. The main construction of two possible bilinear forms on $k\sur B(G,H)$ is performed in Section~6. One form is attached to a given simple $k\Out(H)$-module~$V$, the other corresponds to the largest semi-simple quotient $k\Out(H)/J(k\Out(H))$. We describe the quotient 
$k\sur B(-,H)/R(-,H)$, where $R(G,H)$ is the (right) kernel of the bilinear form. By using the bilinear form corresponding to a given simple $k\Out(H)$-module~$V$, we prove in Section~7 a formula for the dimension of the evaluation at~$G$ of the simple functor $S_{H,V}$. In Section~8, by using the other bilinear form, we prove that $k\sur B(-,H)/R(-,H)$ is the largest semi-simple quotient of the biset functor $k\sur B(-,H)$. The same result may not hold after evaluation at some finite group $G$ and this question is studied in Section~9, where some sufficient conditions are given.
 In Section~10, we analyze the case when the group $H$ is the trivial group, where more information can be obtained. In particular, a rather large subspace of the Jacobson radical $J(k\sur B(G,G))$ is described. In Section~11, we replace the group $H$ by a fixed subquotient $P/Q$ of the group $G$ and describe some natural ideals of $kB(G,G)$ corresponding to this subquotient. The case of the 
subring of $kB(G,G)$ generated by the bisets which are free on both sides is treated in Section~12. Finally, several examples are presented in Section~13.


\section{Biset functors}
\noindent
We first review some know facts about biset functors and refer to~\cite{Bo1} and~\cite{Bo3} for more details. Given two finite groups $G$ and $H$, the Burnside group $B(G,H)$ is the Grothendieck group of the category of finite $(G,H)$-bisets. Since we are interested in simple modules, it is no loss to work over a field~$k$ and we define $kB(G,H)=k\otimes_{\Z}B(G,H)$. In particular, $kB(G,G)$ is a finite dimensional $k$-algebra, called the {\it double Burnside ring\/} of~$G$. We do not need to consider the usual double Burnside ring $B(G,G)$ (defined over~$\Z$) and we work only over the field~$k$. In Section~12, we will consider another version of the double Burnside ring, namely the subring obtained by requiring that the bisets are free on both sides, but for the rest of this paper, we always use the full ring $kB(G,G)$.\par

A {\it section\/} of a finite group $G$ is a pair $(S,T)$ of subgroups of $G$ such that $T$ is a normal subgroup of $S$. In that case, the group $S/T$ is called a {\it subquotient\/} of~$G$. We write $H\sqsubseteq G$ when the group $H$ is isomorphic to a subquotient of~$G$ and we write $H\sqsubset G$ if $H\sqsubseteq G$ and $H\not\cong G$ (hence $|H|<|G|$). We also write $N_G(S,T)$ for the normalizer of the section, that is, the set of all $g\in G$ such that $gSg^{-1}=S$ and $gTg^{-1}=T$.\par

If $(S,T)$ is a section of $G$, then there are elementary bisets $\Res_S^G$, $\Def_{S/T}^S$, $\Ind_S^G$, $\Inf_{S/T}^S$, and their composites $\Defres_{S/T}^G$ and $\Indinf_{S/T}^G$ (see Section~2.3 in~\cite{Bo3}). Also any group isomorphism $\sigma:A\to B$ defines a $(B,A)$-biset $\Iso_\sigma$. 
Recall the following basic result (see Lemma~3 in~\cite{Bo1} or Lemma~2.3.26 in~\cite{Bo3}).

\result{Lemma} \label{biset}
Let $X$ and $Y$ be finite groups.
\begin{enumerate}
\item Any transitive $(X,Y)$-biset has the form $\Indinf_{J/K}^X\, \Iso_{\sigma}\,\Defres_{S/T}^Y$, where $(J,K)$ is a section of $X$, $(S,T)$ is a section of $Y$, and $\sigma:S/T\to J/K$ is a group isomorphism. 
\item Let $\cal E$ be the set of triples $\big( (J,K), \sigma, (S,T) \big)$ where $(J,K)$ is a section of~$X$, $(S,T)$ is a section of $Y$, and $\sigma:S/T\to J/K$ is a group isomorphism. The group $X\times Y$ acts by conjugation on~$\cal E$. Then the set of all elements $\Indinf_{J/K}^X\, \Iso_{\sigma}\,\Defres_{S/T}^Y$ is a $k$-basis of $kB(X,Y)$, where the triple $\big( (J,K), \sigma, (S,T) \big)$ runs over representatives of $(X\times Y)$-orbits in~$\cal E$.
\end{enumerate}
\fresult

The {\it biset category\/} $k\cal C$ is the $k$-linear category whose objects are finite groups, with morphisms $\Hom_{k\cal C}(H,G)=kB(G,H)$ (note that a $(G,H)$-biset is a morphism from $H$ to $G$). The composition of morphisms is the $k$-linear extension of the usual products of bisets $U\times_HV$. Recall that, if $U$ is a $(G,H)$-biset and $V$ is an $(H,L)$-biset, then $U\times_HV$ is a $(G,L)$-biset in the obvious way.\par

A {\it biset functor\/} is a $k$-linear functor from $k\cal C$ to the category $k{-}\Mod$ of $k$-vector spaces and we let $\cal F$ be the category of all such biset functors (an abelian category). We often use a dot for the action of bisets on evaluation of functors, that is, $\alpha\cdot x\in F(G)$ whenever $F\in{\cal F}$, $x\in F(H)$, and $\alpha\in kB(G,H)$. A {\it subquotient\/} of a functor is a quotient of a subfunctor. Moreover, a sequence of functors 
$$0  \;\longrightarrow\; F_1  \;\longrightarrow\; F_2  \;\longrightarrow\; F_3  \;\longrightarrow\; 0$$
is exact if and only if, for every finite group $G$, the evaluation sequence
$$0  \;\longrightarrow\; F_1(G)  \;\longrightarrow\; F_2(G)  \;\longrightarrow\; F_3(G)  \;\longrightarrow\; 0$$
is exact.\par

A biset functor is called {\it simple\/} if it is non-zero and has no proper non-zero subfunctor. The evaluation at a finite group $G$ of a simple functor (and also of a representable functor) is always a finite-dimensional $k$-vector space, so we shall in fact only deal with functors having this additional property.\par

For any fixed finite group $G$, consider the representable functor $kB(-,G)$ (which is a projective functor). Its evaluation at a group $X$ has a natural structure of $(kB(X,X),kB(G,G))$-bimodule. For any $kB(G,G)$-module~$W$, we define, following~\cite{Bo1}, the functor
$$L_{G,W}=kB(-,G)\otimes_{kB(G,G)}W \,,$$
which satisfies the following adjunction property (see Section~2 in~\cite{Bo1}).

\result{Lemma} \label{adjunction1}
Let $G$ be a finite group. The functor
$$kB(G,G){-}\Mod \;\longrightarrow\; {\cal F}  \;,\quad W\mapsto L_{G,W}$$
is left adjoint of the evaluation functor
$${\cal F}  \;\longrightarrow\; kB(G,G){-}\Mod \;,\quad F\mapsto F(G) \,.$$
\fresult

Our next result is a slight extension of the first lemma of~\cite{Bo1}.

\result{Lemma} \label{JGW}
Let $G$ be a finite group, let $W$ be a $kB(G,G)$-module, and let
$$J_{G,W}(X)=\Big\{\sum_i\phi_i\otimes w_i \in L_{G,W}(X) \mid \forall\psi\in kB(G,X), \sum_i (\psi\phi_i)\cdot w_i=0
\Big\} \,.$$
\begin{enumerate}
\item $J_{G,W}$ is the unique subfunctor of $L_{G,W}$ which is maximal with respect to the condition that it vanishes at~$G$.
\item If $W$ is a simple module, then $J_{G,W}$ is the unique maximal subfunctor of $L_{G,W}$ and $L_{G,W}/J_{G,W}$ is a simple functor.
\end{enumerate}
\fresult

\pf
Note first that the condition $\sum_i (\psi\phi_i)\cdot w_i=0$ in~$W$ is equivalent to the condition that $\sum_i\phi_i\otimes w_i$ lies in the kernel of $L_{G,W}(\psi):L_{G,W}(X)\to L_{G,W}(G)$, in view of the isomorphism $L_{G,W}(G)\cong W$ induced by $\phi\otimes w\mapsto \phi w$. Therefore
$$J_{G,W}(X)=\bigcap_{\psi\in kB(G,X)} \Ker\big(L_{G,W}(\psi):L_{G,W}(X)\to L_{G,W}(G) \big) \,.$$
It is then easy to check that $J_{G,W}$ is a subfunctor of~$L_{G,W}$. It is also clear that it vanishes at $G$ since $J_{G,W}(G)$ lies in the kernel of~$\Id$. Finally, if $F$ is a subfunctor of $L_{G,W}$ vanishing at~$G$, then $F(X)$ must lie in the kernel of all maps $L_{G,W}(\psi)$ for $\psi\in kB(G,X)$, so that $F\subseteq J_{G,W}$.\par

The case where $W$ is simple is explicit in~\cite{Bo1} and we just recall the main line of the argument. Any subfunctor $F$ of $L_{G,W}$ either vanishes at~$G$, hence $F\subseteq J_{G,W}$, or satisfies $F(G)=W$, hence $F=L_{G,W}$. Therefore $J_{G,W}$ is the unique maximal subfunctor of $L_{G,W}$ and $L_{G,W}/J_{G,W}$ is a simple functor.
\endpf

We define $S_{G,W}=L_{G,W}/J_{G,W}$ and we emphasize that this is a simple functor provided $W$ is a simple $kB(G,G)$-module. This provides the first link between simple $kB(G,G)$-modules and simple functors. Our next section shows that the connection is much stronger.


\section{Simple modules and simple functors}
\noindent
We prove in this section that the link between simple $kB(G,G)$-modules and simple biset functors is deep.
Note that the results of this section work in the more general context presented in Section~3.2 of~\cite{Bo3},  in particular for inflation functors (corresponding to bisets which are free on one side), for global Mackey functors (corresponding to bisets which are free on both sides), and also for functors defined only on a given class of finite groups. For simplicity, we consider all groups and all bisets. An overview of the case of global Mackey functors, using bisets which are free on both sides, is given in Section~12.\par

Our first proposition is a special case of the results of Section~4.2 in~\cite{Bo3}. We repeat the arguments for convenience.

\result{Proposition} \label{simple-evaluation}
\begin{enumerate}
\item If $G$ is a finite group and $W$ is a simple $kB(G,G)$-module, then $W$ is the evaluation at~$G$ of a simple biset functor, namely $S_{G,W}$. Moreover $S_{G,W}$ is the unique simple functor, up to isomorphism,  such that its evaluation at~$G$ is isomorphic to~$W$ as a $kB(G,G)$-module.
\item If $S$ is a simple biset functor and $G$ is a finite group, then either $S(G)=0$ or $S(G)$ is a simple $kB(G,G)$-module. In the latter case, $S\cong S_{G,W}$, where $W=S(G)$.
\end{enumerate}
\fresult

\pf
The first claim of (1) is an immediate consequence of~Lemma~\ref{JGW} above, because $S_{G,W}(G)=L_{G,W}(G)\cong W$. The second claim of (1) is a consequence of (2).\par

For the proof of (2), suppose that $S(G)\neq 0$. Let $M$ be a non-zero $kB(G,G)$-submodule of~$S(G)$ and let $i: M\to S(G)$ be the inclusion. By the adjunction of Lemma~\ref{adjunction1}, this corresponds to a non-zero morphism $\theta:L_{G,M}\to S$. Since $S$ is simple, $\theta$ must be surjective, that is, surjective on every evaluation. Thus $\theta(G): L_{G,M}(G)\to S(G)$ is surjective, but this is just $\theta(G)=i: M\to S(G)$. This proves that $M=S(G)$. Consequently $S(G)$ is simple. Letting now $W=S(G)$, we obtain a surjective morphism $\theta:L_{G,W}\to S$. But $S_{G,W}$ is the unique simple quotient of $L_{G,W}$, by Lemma~\ref{JGW} above, so $S\cong S_{G,W}$.
\endpf

The fact that a simple $kB(G,G)$-module determines uniquely a simple biset functor shows that the study of simple $kB(G,G)$-modules is, in some sense, equivalent to the study of simple biset functors. So we have, in some sense, enriched the structure of objects we are working with. On the other hand, the evaluations of simple functors are unfortunately not easy to determine. One of our main results in Section~7 will actually give some answer to this question.\par

A simple functor $S$ has many realizations $S\cong S_{G,W}$, one for each $G$ such that $S(G)\neq 0$. But recall that, in order to obtain a parametrization, one can do better, as follows (see Section~4 in~\cite{Bo1} or Section~4.3 in~\cite{Bo3}).

\result{Proposition} \label{parametrization}
Let $S$ be a simple biset functor, let $H$ be a group of minimal order such that $S(H)\neq 0$, and let $V=S(H)$.
\begin{enumerate}
\item $H$ is unique up to isomorphism and $S\cong S_{H,V}$.
\item Let $kI(H,H)$ be the ideal of $kB(H,H)$ generated by all bisets which factor through a proper subquotient of~$H$, so that $kB(H,H)/kI(H,H)\cong k\Out(H)$ (where $\Out(H)$ denotes the group of outer automorphisms of~$H$). Then $kI(H,H)$ acts by zero on~$V$ and $V$ is a $k\Out(H)$-module.
\item If $S(G)\neq0$ for some finite group $G$, then $H\sqsubseteq G$. 
\end{enumerate}
\fresult

This provides a {\it parametrization\/} of simple functors by pairs $(H,V)$ where $H$ is a finite group and $V$ is a simple $k\Out(H)$-module. In the next sections of this paper, we shall often fix a finite group $H$ and a simple $k\Out(H)$-module~$V$. Then we shall work with the simple functor $S_{H,V}$ and consider finite groups $G$ such that $W=S_{H,V}(G)\neq0$ (so that in fact $S_{H,V}\cong S_{G,W}$).\par

\result{Corollary} \label{counting} 
If $G$ is a finite group, the number of isomorphism classes of simple $kB(G,G)$-modules is equal to the number of isomorphism classes of pairs $(H,V)$, where $H$ is finite group such that $H\sqsubseteq G$ and $V$ is a simple $k\Out(H)$-module, subject to the condition that $S_{H,V}(G)\neq0$.
\fresult

\pf
This follows immediately from Propositions~\ref{simple-evaluation} and \ref{parametrization}.
\endpf

It follows from this proposition that the question of the vanishing of evaluation of simple functors is a crucial issue for the description of simple $kB(G,G)$-modules.  Whenever $S_{H,V}(G)=0$, the pair $(H,V)$ must not be counted. This question of the vanishing of evaluation of simple functors will be considered in another paper~\cite{BST}.

Another useful fact is the following~:

\result{Proposition} \label{generated}
Let $S$ be a simple biset functor and let $G$ be a group such that $S(G)\neq 0$. Then $S$ is generated by $S(G)$, that is, $S(X)=kB(X,G)S(G)$ for all finite groups~$X$. More precisely, if $0\neq u\in S(G)$, then $S(X)=kB(X,G)\cdot u$.
\fresult

\pf
Given $0\neq u\in S(G)$, let $S'(X)=kB(X,G)\cdot u$ for all finite groups~$X$. This clearly defines a non-zero subfunctor $S'$ of~$S$, so $S'=S$ by simplicity of~$S$.
\endpf

The connection between simple functors and simple evaluations extends further, as follows.

\result{Proposition} \label{subquotients}
Let $S$ be a simple biset functor and let $G$ be a finite group such that $S(G)\neq 0$. Let $F$ be any biset functor. Then the following are equivalent:
\begin{enumerate}
\item $S$ is isomorphic to a subquotient of~$F$.
\item The simple $kB(G,G)$-module $S(G)$ is isomorphic to a subquotient of the $kB(G,G)$-module $F(G)$.
\end{enumerate}
\fresult

\pf
It is clear that (1) implies (2). Suppose that (2) holds and let $W_1$, $W_2$ be submodules of~$F(G)$ such that $W_2\subset W_1$ and $W_1/W_2\cong S(G)$. For $i\in\{1,2\}$, let $F_i$ be the subfunctor of $F$ generated by $W_i$. Explicitly, for any finite group $X$, 
$F_i(X)=kB(X,G)\cdot W_i\subseteq F(X)$. Then $F_i(G)=W_i$ and $(F_1/F_2)(G)=W_1/W_2\cong S(G)$. The isomorphism $S(G)\to (F_1/F_2)(G)$ induces, by the adjunction of Lemma~\ref{adjunction1}, a non-zero morphism $\theta: L_{G,S(G)}\to F_1/F_2$. Since $S(G)$ is simple, $L_{G,S(G)}$ has a unique maximal subfunctor $J_{G,S(G)}$, by Lemma~\ref{JGW}, and $L_{G,S(G)}/J_{G,S(G)} \cong S_{G,S(G)} = S$, by Proposition~\ref{simple-evaluation}. Let $F'_1=\theta(L_{G,S(G)})$ and $F'_2=\theta(J_{G,S(G)})$. Since $\theta\neq0$, we obtain
$$F'_1/F'_2 \cong L_{G,S(G)}/J_{G,S(G)} \cong S_{G,S(G)} = S \,,$$
showing that $S$ is isomorphic to a subquotient of~$F$.
\endpf

Note that, in the proof above, we have $F_2\subseteq F'_2\subseteq F'_1\subseteq F_1\subseteq F$.
Moreover $(F'_1/F'_2)(G)\cong S(G)\cong (F_1/F_2)(G)$, so that we have $(F_1/F'_1)(G)=0$ and ${(F'_2/F_2)(G)=0}$. But in general we don't have $F'_i=F_i$.


\section{The standard quotient of a representable functor}
\noindent
In this section, we construct the biset functor $k\sur B(-,H)$, which plays a central role in the rest of the paper.
Let us fix a finite group $H$ and consider the representable functor $kB(-,H)$. For every finite group $X$, define
$$I(X,H):=\sum_{K\sqsubset H}B(X,K)B(K,H) \,.$$
This is a subgroup of the abelian group $B(X,H)$. Extending scalars to $k$, we obtain the $k$-subspace
$$kI(X,H):=\sum_{K\sqsubset H}kB(X,K)B(K,H) \,.$$

\result{Lemma} \label{subfunctor}
$kI(-,H)$ is a subfunctor of $kB(-,H)$.
\fresult

\pf
We have
\begin{eqnarray*}
B(Y,X)I(X,H)\;=&B(Y,X)\sum_{K\sqsubset H}B(X,K)B(K,H)\\
=&\sum_{K\sqsubset H}B(Y,X)B(X,K)B(K,H) \\
\subseteq&\sum_{K\sqsubset H}B(Y,K)B(K,H) =I(Y,H) \,.
\end{eqnarray*}
\endpf

The idea of passing to the quotient by all morphisms factorizing below~$H$ has been widely used (e.g. in Section~4 of~\cite{We2}).
We define the {\it standard quotient\/} of $kB(-,H)$ to be
$$k\sur B(-,H) =kB(-,H)/kI(-,H)\,.$$
This terminology will be motivated below.
In particular $kI(H,H)$ is an ideal of the double Burnside ring $kB(H,H)$ and $k\sur B(H,H)\cong k\Out(H)$, as already noticed in Proposition~\ref{parametrization}.\par

We first note the following elementary result.

\result{Lemma} \label{sqsubset}
Let $H$ and $X$ be finite groups. Then $\sur B(X,H) \neq 0$ if and only if $H\sqsubseteq  X$.
\fresult

\pf
If $\sur B(X,H) \neq 0$, there exists a transitive $(X,H)$-biset $U$ such that its class in $\sur B(X,H)$ is non-zero. But by Lemma~\ref{biset}, we can write
$$U=\Indinf_{J/K}^X \Iso_{\sigma}\Defres_{S/T}^H \,,$$
where $(S,T)$ is a section of $H$, $(J,K)$ is a section of $X$, and $\sigma:S/T\to J/K$ is an isomorphism. If $S/T$ was a proper subquotient of~$H$, then we would have $\Defres_{S/T}^H\in I(S/T,H)$, hence $U\in I(X,H)$, and the class of~$U$ would be zero in $\sur B(X,H)$. Therefore $S=H$ and $T=1$, so that $U=\Indinf_{J/K}^X \Iso_{\sigma}$ where $\sigma:H\to J/K$ is an isomorphism. This proves that $H\sqsubseteq  X$.\par

Conversely, if $H\sqsubseteq  X$, then there is an isomorphism $\sigma:H\to J/K$ where $(J,K)$ is a section of $X$, and $U=\Indinf_{J/K}^X \Iso_{\sigma}$ defines a non-zero element in $\sur B(X,H)$, by Lemma~\ref{biset}.
\endpf

For any finite group $X$, the evaluation $k\sur B(X,H)$ has a natural structure of right $k\Out(H)$-module, because the right action of $kI(H,H)$ is zero. For any left $k\Out(H)$-module~$V$, we define the functor
$$\sur L_{H,V}=k\sur B(-,H)\otimes_{k\Out(H)}V \,,$$
or also $\sur L_{H,V}\cong k\sur B(-,H)\otimes_{kB(H,H)}V$
if $V$ is viewed as a $kB(H,H)$-module. Applying the functor $-\otimes_{kB(H,H)}V$ to the exact sequence
$$0\to kI(X,H) \to kB(X,H) \to k\sur B(X,H) \to 0 \,,$$
we see that $\sur L_{H,V}(X)$ is a quotient of $L_{H,V}(X)$ (quotient by the image of $kI(X,H)\otimes_{kB(H,H)}V$). We call $\sur L_{H,V}$ the {\it standard quotient\/} of $L_{H,V}$. Note that if $V=k\Out(H)$, then $\sur L_{H,k\Out(H)}=k\sur B(-,H)$.\par

Following Section~4 of~\cite{We2}, we define, for any biset functor $F$ and any finite group $H$, the {\it restriction kernel\/}
$$\sou F(H)=\bigcap_{K\sqsubset H,\,\phi\in B(K,H)} \Ker(F(\phi)) \,.$$
Clearly $kI(H,H)$ acts by zero on
$\sou F(H)$, so that $\sou F(H)$ is a $k\Out(H)$-module. We also refer to~\cite{Ya} for a use of restriction kernels. Similarly to the adjunction of Lemma~\ref{adjunction1}, we have the following analogous property for the standard quotient.

\result{Lemma} \label{adjunction2}
Let $H$ be a finite group. The functor
$$k\Out(H){-}\Mod \;\longrightarrow\; {\cal F}  \;,\quad V\mapsto \sur L_{H,V}$$
is left adjoint of the functor
$${\cal F}  \;\longrightarrow\; k\Out(H){-}\Mod \;,\quad F\mapsto \sou F(H) \,.$$
\fresult

\pf
Let $g: V\to \sou F(H)$ be a $k\Out(H)$-linear map. Composing with the inclusion $i:\sou F(H)\to F(H)$ and applying the adjunction of Lemma~\ref{adjunction1}, we obtain a morphism $\theta:L_{H,V}\to F$ such that $\theta(H)=i\circ g$. Now fix $K\sqsubset H$ and let $\alpha\in B(X,K)$ and $\beta\in B(K,H)$, so that $\alpha\beta$ is among the generators of $kI(X,H)$. Then we obtain
$$\theta(X)L_{H,V}(\alpha\beta) = F(\alpha\beta)\theta(H)=F(\alpha)F(\beta)i g=0 \,,$$
because $F(\beta)i=0$ by definition of $\sou F(H)$. It follows that $\theta(X)$ is zero on the image of $kI(X,H)\otimes_{kB(H,H)}V$ in $L_{H,V}(X)=kB(X,H)\otimes_{kB(H,H)}V$. Therefore $\theta(X)$ induces a morphism $\sur\theta(X):\sur L_{H,V}(X) \to F(X)$. Then $\sur\theta:\sur L_{H,V} \to F$ is the morphism corresponding to~$g$ under the required adjunction.
\endpf

In Section~5 of~\cite{We2}, Webb defines a biset functor $\Delta_{H,V}$ and proves that $V\mapsto \Delta_{H,V}$ is left adjoint of the functor $F\mapsto \sou F(H)$. Therefore we have in fact
$$\Delta_{H,V} \cong \sur L_{H,V} \,.$$
One of the main results of Webb asserts that the functors $\Delta_{H,V}$, where $V$ is simple, are the {\it standard objects\/} in a highest weight category structure on~$\cal F$ (over a field of characteristic zero). Thus we see that the functors $\sur L_{H,V}$, hence also the functors $k\sur B(-,H)$, deserve to be called standard. This explains our terminology of {\it standard quotient\/}.\par

Webb also shows that, if $V$ is a simple $k\Out(H)$-module, then the simple functor $S_{H,V}$ is a quotient of $\Delta_{H,V}$. We now show that this works also for any $k\Out(H)$-module $V$, using our definition of the corresponding functor $S_{H,V}$.

\result{Proposition} \label{simple-quotient}
Let $H$ be a finite group and $V$ a $k\Out(H)$-module.
\begin{enumerate}
\item The image of $kI(X,H)\otimes_{kB(H,H)}V$ in $kB(X,H)\otimes_{kB(H,H)}V=L_{H,V}(X)$ is contained in $J_{H,V}(X)$.
\item The quotient morphism $\varepsilon:kB(-,H)\to k\sur B(-,H)$ induces a surjective morphism $\varepsilon:L_{H,V}\to \sur L_{H,V}$ and an isomorphism
$$\sur\varepsilon:S_{H,V}=L_{H,V}/J_{H,V} \flh{\sim}{} \sur L_{H,V}/\sur J_{H,V} \,,$$
where $\sur J_{H,V}=\varepsilon(J_{H,V})$.
\item If $V$ is a simple module, then $\sur J_{H,V}$ is the unique maximal subfunctor of $\sur L_{H,V}$.
\end{enumerate}
\fresult

\pf
Let $\alpha\in kI(X,H)$. For any $\psi\in kB(H,X)$, we have $\psi\alpha\in kI(H,H)$, hence $\psi\alpha\cdot v=0$ for any $v\in V$. This shows that $\alpha\otimes v\in J_{H,V}(X)$, proving the first part. The other parts follow immediately.
\endpf

\result{Remark} \label{surJ} {\rm
By copying the proof of Lemma~\ref{JGW}, one easily obtains that
$$\sur J_{H,V}(X)=\Big\{\sum_i\sur\phi_i\otimes v_i \in \sur L_{G,W}(X) \mid \forall\psi\in kB(H,X), \sum_i (\sur{\psi\phi}_i)v_i=0 \Big\} \,,$$
and that $\sur J_{H,V}$ is the unique subfunctor of $\sur L_{H,V}$ which is maximal with respect to the condition that it vanishes at~$H$.
}
\fresult


\section{The minimal group of a simple $kB(G,G)$-module}
\noindent
In this section, we fix a finite group $G$ and we let $W$ be a simple $kB(G,G)$-module. A {\it minimal group\/} for $W$ is a group $P$ of minimal order subject to the condition that $B(G,P)B(P,G)W\neq0$. Note that composition of bisets defines a homomorphism $B(G,P)B(P,G)\to B(G,G)$, so that $B(G,P)B(P,G)$ acts on~$W$. Note also that the image of this homomorphism is a two-sided ideal of $B(G,G)$, generated by all the $(G,G)$-bisets which factor through~$P$, so that $B(G,P)B(P,G)W$ is a $kB(G,G)$-submodule of $W$ (hence either $0$ or $W$ by simplicity of~$W$). This notion of minimality only involves the double Burnside ring, but it has a very useful interpretation in terms of biset functors, as follows.

\result{Proposition} \label{minimal}
Let $G$ be a finite group, let $W$ be a simple $kB(G,G)$-module, and let $P$ be a minimal group for~$W$.
\begin{enumerate}
\item$P$ is unique up to isomorphism and $P\sqsubseteq G$. 
\item Let $(H,V)$ be a pair which parametrizes the simple functor $S_{G,W}$, that is, $S_{G,W}=S_{H,V}$, where $H$ is a group of minimal order such that ${S_{G,W}(H)\neq0}$ and $V$ is a simple $k\Out(H)$-module (see Proposition~\ref{parametrization}). Then $H$ is a minimal group for~$W$. Moreover $V=S_{G,W}(H)$ and $W=S_{G,W}(G)$.
\end{enumerate}
\fresult

\pf
We first give a proof using the simple functor $S=S_{G,W}=S_{H,V}$,  where $(H,V)$ is as in statement~(2), so in particular $H\sqsubseteq G$. We have $W=S(G)$ and $B(P,G)W=B(P,G)S(G)$ makes sense because bisets can be applied to this evaluation. If $B(G,P)B(P,G)W\neq0$, then $B(P,G)S(G)\neq0$, hence $S(P)\neq 0$. This implies that $H$ is isomorphic to a subquotient of~$P$, by Proposition~\ref{parametrization}. On the other hand, by Proposition~\ref{generated}, we have
$$B(G,H)B(H,G)W=B(G,H)B(H,G)S(G)=B(G,H)S(H)=S(G)\neq0\,.$$
By minimality of $P$, this implies that the subquotient $H$ of~$P$ must be isomorphic to~$P$. This proves both statements.\par

We now sketch a second proof of the first statement, using only bisets. 
Let $\alpha\in B(G,P)$ and $\beta\in B(P,G)$ be transitive bisets such that $\alpha\beta W\neq0$. By minimality of~$P$, $\alpha\beta$ cannot factor through a proper subquotient of~$P$. Therefore, by Lemma~\ref{biset}, $\alpha=\Indinf_{J/K}^G \Iso_{\sigma}$ for some section $(J,K)$ of~$G$ and $\beta=\Iso_{\rho} \Defres_{S/T}^G$ for some section $(S,T)$ of~$G$.
If now $P'$ is another minimal group for~$W$, then
$$B(G,P')B(P',G)B(G,P)B(P,G)W=W$$
and we get transitive bisets $\alpha'\in B(G,P')$ and $\beta'\in B(P',G)$ such that
$$\alpha'\beta' \alpha\beta W\neq0 \,.$$
Then we decompose the $(P',P)$-biset $\beta' \alpha$ as a sum of transitive bisets and we obtain at least one transitive summand $\gamma\in B(P',P)$ such that $\alpha'\gamma\beta W\neq0$. But $\gamma$ factorizes through subquotients of $P$ and $P'$ (by Lemma~\ref{biset}), so by minimality of $P$ and $P'$, the only possibility is that $\gamma$ is an isomorphism. Therefore $P\cong P'$, as required.
\endpf
 
We now establish the link bewteen minimal groups and standard quotients.

\result{Proposition} \label{quotient}
Let $G$ be a finite group, let $W$ be a simple $kB(G,G)$-module, and let $P$ be a minimal group for~$W$.
\begin{enumerate}
\item $W$ is isomorphic to a quotient of the $kB(G,G)$-module $k\sur B(G,P)$.
\item If $W$ is isomorphic to a quotient of the $kB(G,G)$-module $k\sur B(G,H)$ for some finite group $H$, then $H$ is isomorphic to a subquotient of~$P$.
\end{enumerate}
\fresult

\pf
By Proposition~\ref{minimal}, $P$ is a minimal group for $S_{G,W}$ and we have $S_{G,W}=S_{P,V}$, where $V$ is a simple $k\Out(P)$-module. Moreover, there are surjective morphisms
$$k\sur B(-,P) \longrightarrow \sur L_{P,V} \longrightarrow S_{P,V} \,,$$
where the first morphism maps $\phi\in k\sur B(X,P)$ to $\phi\otimes v\in\sur L_{P,V}(X)$ for some fixed $v\in V$, and the second comes from Proposition~\ref{simple-quotient}. By evaluating at~$G$, we see that $W=S_{G,W}(G)=S_{P,V}(G)$ is isomorphic to a quotient of the $kB(G,G)$-module $k\sur B(G,P)$, proving~(1).\par

Assume now that $W=S_{P,V}(G)$ is isomorphic to a quotient of the $kB(G,G)$-module $k\sur B(G,H)$ for some~$H$.
By Proposition~\ref{subquotients}, $S_{P,V}$ is isomorphic to a subquotient of $k\sur B(-,H)$. By evaluating at~$P$, we see that $V=S_{P,V}(P)$ is isomorphic to a subquotient of the $kB(G,G)$-module $k\sur B(P,H)$.
Therefore $k\sur B(P,H)\neq0$, hence $H\sqsubseteq P$ by  Proposition~\ref{sqsubset}, proving~(2).
\endpf

The above results show that it is worth considering the simple quotients of $k\sur B(G,H)$ in order to find the simple $kB(G,G)$-modules. Every such simple module must appear, possibly for several finite groups~$H$. Our next task is to study further semi-simple quotients of $k\sur B(G,H)$.


\section{Bilinear forms on standard quotients}
\noindent
In this section, we introduce one of the main constructions of this paper.
We fix a finite group $H$ and consider a quotient algebra $E$ of $k\Out(H)$ with corresponding $k$-algebra map $\sur\pi: k\Out(H) \to E$, which we compose with the quotient map $kB(H,H)\to k\Out(H)$ to get
$$\pi : kB(H,H) \longrightarrow E \,.$$
We assume that $E$ is a symmetric algebra with symmetrizing form $\tau:E\to k$. This means that $\tau$ is $k$-linear and induces a symmetric non-degenerate bilinear form on~$E$ given by
$$(a,b)\mapsto \tau(ab) \,.$$
We have in mind two cases which we shall consider later~:

\begin{enumerate}
\item $k$ is algebraically closed (or at least large enough for the group $\Out(H)$), $E=\End_k(V)$ where $V$ is a simple $k\Out(H)$-module, $\sur\pi: k\Out(H)\to \End_k(V)$ is the natural surjection, and $\tau=\tau_V$ is the trace map on $\End_k(V)$. Note that if $V$ is self-dual, then $\tau_V$ satisfies the additional condition $\tau_V(\sur\pi(s^{-1}))= \tau_V(\sur\pi(s))$ for all $s\in\Out(H)$.
\item $k$ is algebraically closed (or at least large enough for the group $\Out(H)$), $E=k\Out(H)/J(k\Out(H))\cong \prod_{i=1}^r \End_k(V_i)$, where $V_1,\ldots,V_r$ are the simple $k\Out(H)$-modules, $\sur\pi: k\Out(H)\to k\Out(H)/J(k\Out(H))$ is the quotient map, and $\tau=\sum_{i=1}^r \tau_{V_i}$, the sum of the trace maps on each $\End_k(V_i)$. Note that $\tau$ satisfies the additional condition $\tau(\sur\pi(s^{-1}))= \tau(\sur\pi(s))$ for all $s\in\Out(H)$.
If the characteristic of $k$ does not divide $|\Out(H)|$, then $E=k\Out(H)$, $\sur\pi=\Id$, and we could also choose for $\tau$ the ordinary map for group algebras (coefficient of~1), because any symmetrizing form would do.
\end{enumerate}

With such data, we can construct a bilinear form on $k\sur B(X,H)$, for any finite group $X$.
Recall that any $(X,H)$-biset $U$ has an opposite $U\op$, which is an $(H,X)$-biset (see Section~2.3 in~\cite{Bo3}), and this extends to a $k$-linear map $kB(X,H)\to kB(H,X)$, $\alpha\mapsto\alpha\op$, such that $(\alpha\beta)\op=\beta\op\alpha\op$. Then we define
$${<}-,-{>}_X \,: k\sur B(X,H) \times k\sur B(X,H) \longrightarrow k \;,
\qquad {<}\sur\phi,\sur\psi{>}_X = \tau \pi(\phi\op \psi) \,,$$
where $\phi,\psi\in kB(X,H)$ are representatives of $\sur\phi$ and $\sur\psi$ respectively.

\result{Lemma} \label{bilinear}
Let $H$ and $X$ be finite groups and let $\pi: kB(H,H)\to E$ and $\tau:E\to k$ be as above. Let $R(X,H)$ be the right kernel of the form ${<}-,-{>}_X$, that is, the set of all elements $\sur\psi\in k\sur B(X,H)$ such that  ${<}\sur\phi,\sur\psi{>}_X =0$ for all $\sur\phi\in k\sur B(X,H)$,
\begin{enumerate}
\item The map ${<}-,-{>}_X$ is well-defined and $k$-bilinear.
\item If $\tau$ satisfies the condition $\tau(\sur\pi(s^{-1}))= \tau(\sur\pi(s))$ for all $s\in\Out(H)$ (e.g. in case~(2) above), then the bilinear form ${<}-,-{>}_X$ is symmetric.
\item We have ${<}\sur{\alpha\op}\,\sur\beta,\sur\psi{>}_X = {<}\sur\beta,\sur\alpha\,\sur\psi{>}_Y$ for all $\sur\alpha\in k\sur B(Y,X)$, $\sur\beta\in k\sur B(Y,H)$, $\sur \psi\in k\sur B(X,H)$.
\item $R(-,H)$ is a subfunctor of $k\sur B(-,H)$.
\item We have ${<}\sur\phi,\sur\psi\,\sur\gamma{>}_X ={<}\sur\phi\,\sur\gamma\op,\sur\psi{>}_X$ for all $\sur\phi,\sur\psi\in k\sur B(X,H)$ and $\sur\gamma\in k\sur B(H,H)=k\Out(H)$.
\item $R(X,H)$ is a right $k\Out(H)$-submodule of $k\sur B(X,H)$.
\item $R(H,H) = \Ker(\sur\pi)$.
\end{enumerate}
\fresult

\pf
(1) Replace $\psi\in kB(X,H)$ by $\psi'=\psi+\chi$, where $\chi\in kI(X,H)$. Then, for all $\phi\in kB(X,H)$,
 $$\pi(\phi\op\psi')=\pi(\phi\op\psi)+\pi(\phi\op\chi)=\pi(\phi\op\psi) \,,$$
because $\phi\op\chi\in kI(H,H)$, hence $\pi(\phi\op\chi)=0$. A similar argument applies if we modify $\phi$ by an element of $kI(X,H)$ and this shows that the form is well-defined. It is obvious that it is bilinear.\par

(2) Let $\delta=\phi\op \psi\in kB(H,H)$, so that $\psi\op \phi=\delta\op$. The map $\pi$ factors through $k\Out(H)$ and the passage to opposite bisets induces on $k\Out(H)$ the map $s\mapsto s^{-1}$ for each $s\in\Out(H)$. Since we have $\tau(\sur\pi(s^{-1}))= \tau(\sur\pi(s))$ by our assumption, we obtain $\tau\pi(\delta)=\tau\pi(\delta\op)$, and the symmetry follows.\par

(3) This follows immediately from the equality $(\alpha\op\beta)\op\psi=\beta\op(\alpha\psi)$.\par

(4) Let $\sur\psi\in R(X,H)$, so that ${<}\sur\phi,\sur\psi{>}_X=0$ for all $\sur\phi\in k\sur B(X,H)$. Apply the morphism $\alpha\in kB(Y,X)$.
Then, for all $\sur\beta\in k\sur B(Y,H)$, we have
${<}\sur\beta,\sur\alpha\,\sur\psi{>}_Y = {<}\sur{\alpha\op}\,\sur\beta,\sur\psi{>}_X=0$.
Therefore $\alpha\cdot \sur\psi\in R(Y,H)$. This shows that $R(-,H)$ is a subfunctor of $k\sur B(-,H)$.\par

(5) Using the symmetry of the form~$\tau$, we obtain 
$${<}\sur\phi,\sur\psi\,\sur\gamma{>}_X =\tau\pi(\phi\op\psi\gamma)=\tau\pi(\gamma\phi\op\psi)
= \tau\pi((\phi\gamma\op)\op\psi)={<}\sur\phi\,\sur\gamma\op,\sur\psi{>}_X \,.$$

(6) Let $\sur\psi\in R(X,H)$, and let $\sur\gamma\in k\sur B(H,H)=k\Out(H)$. Then, for all $\sur\phi\in k\sur B(X,H)$,
$${<}\sur\phi,\sur\psi\,\sur\gamma{>}_X={<}\sur\phi\,\sur\gamma\op,\sur\psi{>}_X =0$$
This proves that $\sur\psi\,\sur\gamma\in R(X,H)$.

(7) Let $\sur\gamma\in k\sur B(H,H)=k\Out(H)$. Then $\sur\gamma\in R(H,H)$ if and only if $\tau \pi(\phi\op \gamma)=0$ for all $\phi\in kB(H,H)$. But this gives $\tau \big(\pi(\phi\op) \pi(\gamma)\big)=0$ for all~$\phi$, hence $\pi(\gamma)=0$ by non-degeneracy of the bilinear form on~$E$ induced by~$\tau$. Therefore $\sur\gamma\in R(H,H)$ if and only if $\sur\gamma\in \Ker(\sur\pi)$.
\endpf

If necessary, it could be possible to define the bilinear form directly on $kB(X,H)$ instead of its quotient $k\sur B(X,H)$. Then $kI(X,H)$ would be in the kernel of the form and the form would induce the one defined above.\par

We are interested in the quotient $k\sur B(-,H)/R(-,H)$ and we first determine its structure.

\result{Theorem} \label{B/R}
Let $H$ be a finite group, let $\pi: kB(H,H)\to E$ and $\tau:E\to k$ be as above. Let $R(-,H)$ be the subfunctor of $k\sur B(-,H)$ defined in Lemma~\ref{bilinear}.
Then the quotient functor $k\sur B(-,H)/R(-,H)$ is isomorphic to $S_{H,E}$.
\fresult

\pf
Consider the exact sequence $0\to \Ker(\sur\pi) \to k\Out(H) \to E \to 0$.
Tensor this sequence with $k\sur B(X,H)$, where $X$ is a finite group.
Since tensoring is right exact and since $k\sur B(X,H)\otimes_{k\Out(H)} k\Out(H) \cong k\sur B(X,H)$, we obtain
$$0 \longrightarrow k\sur B(X,H)\cdot \Ker(\sur\pi) \longrightarrow k\sur B(X,H) \longrightarrow
k\sur B(X,H)\otimes_{k\Out(H)} E \longrightarrow 0 \,,$$
where the first map is just the inclusion map.
Note that $k\sur B(X,H)\otimes_{k\Out(H)} E = \sur L_{H,E}(X)$.
Now it is obvious that $k\sur B(X,H)\cdot \Ker(\sur\pi)$ is contained in the kernel $R(X,H)$ of the bilinear form. Therefore we obtain the exact sequence
$$0 \longrightarrow k\sur B(X,H)\cdot \Ker(\sur\pi) \longrightarrow R(X,H) \longrightarrow
R(X,H)\otimes_{k\Out(H)} E \longrightarrow 0 \,.$$
We claim that $R(X,H)\otimes_{k\Out(H)} E$ is equal to $\sur J_{H,E}(X)$,
as defined in Proposition~\ref{simple-quotient}. It then follows that
\begin{eqnarray*}
k\sur B(X,H)/R(X,H) &\cong&
\big(k\sur B(X,H)\otimes_{k\Out(H)} E\big) \big/ \big(R(X,H)\otimes_{k\Out(H)} E\big) \\
&=& \sur L_{H,E}(X)/\sur J_{H,E}(X) \\
& \cong& S_{H,E}(X) \,,
\end{eqnarray*}
by Proposition~\ref{simple-quotient}, proving the first statement.\par

Now we prove the claim. Any element of $k\sur B(X,H)\otimes_{k\Out(H)} E$ can be written $\sur\phi\otimes 1_E$ with $\phi \in kB(X,H)$, because $\sur\pi: k\Out(H) \to E$ is surjective. Then we have

\begin{eqnarray*}
\sur\phi\in R(X,H) \; &\Leftrightarrow& \; {<}\sur\alpha,\sur\phi{>}_X=0 \,, \forall \alpha\in kB(X,H) \\
&\Leftrightarrow&\; \tau\pi(\beta\phi)=0 \,, \forall \beta\in kB(H,X) \\
&\Leftrightarrow&\; \tau\big(\pi(\gamma)\pi(\beta\phi)\big)=0 \,, \forall \beta\in kB(H,X) \,,\forall \gamma\in kB(H,H) \\
&\Leftrightarrow&\; \pi(\beta\phi)=0 \,, \forall \beta\in kB(H,X) \\
&\Leftrightarrow&\; \sur{\beta\phi}\cdot 1_E=0 \,, \forall \beta\in kB(H,X) \\
&\Leftrightarrow&\; \sur\phi\otimes 1_E\in \sur J_{H,E}(X) \,,
\end{eqnarray*}
where we used the non-degeneracy of the bilinear form on~$E$ induced by~$\tau$.
This completes the proof.
\endpf

\result{Corollary} \label{E-semi-simple}
With the notation of Theorem~\ref{B/R}, if $E$ is semi-simple as a $k\Out(H)$-module, then the biset functor $k\sur B(-,H)/R(-,H)$ is semi-simple.
\fresult

\pf
This follows from the observation that $S_{H,V\oplus W}\cong S_{H,V}\oplus S_{H,W}$ and that $S_{H,V}$ is a simple functor if $V$ is a simple $k\Out(H)$-module.
\endpf


\section{Evaluation of simple functors}
\noindent
In this section, we use the results of Section~6 to determine the dimension of the evaluation of a simple functor $S_{H,V}$. In general, the explicit computation of an evaluation $S_{H,V}(G)$ is not easy. There are general procedures for the determination of such evaluations (see Theorem~4.3.20 in~\cite{Bo3} and also  Section~11.2 in~\cite{We2}) and there are some special cases where the evaluation is known in detail, for instance when $k$ has characteristic zero, $H$ is a $b$-group, and $V=k$ is the trivial module (see Section~7.2.4 in~\cite{Bo1} and Theorem~5.5.4 in~\cite{Bo3}).\par

Here we prove that the dimension of an evaluation $S_{H,V}(G)$ can be obtained as the rank of a bilinear form, a result which had been known for a long time in the case where $V=k$, the trivial $k\Out(H)$-module (see Section~8.2 in~\cite{Bo1}), but remained open in the general case.
In Section~11.2 of~\cite{We2}, Webb mentions a bilinear form for computing $\dim(S_{H,V}(G))$ (although this is made explicit only if $\dim(V)=1$), but the form is defined on another vector space. However, the result sketched by Webb is quite similar to our theorem below.\par

The construction of the simple functor $S_{H,V}$ uses a quotient of $L_{H,V}$ or $\sur L_{H,V}$, hence one of the tensor products
$$L_{H,V} =kB(-,H)\otimes_{kB(H,H)} V \quad\text{or}\quad
\sur L_{H,V} =k\sur B(-,H)\otimes_{k\Out(H)} V \,.$$
The advantage of Theorem~\ref{B/R} is that we obtain $S_{H,E}$ directly as a quotient of $k\sur B(-,H)$, without needing to tensor with~$E$. This implies that the dimension of $S_{H,E}(G)$ can be computed by a direct use of the bilinear form and its kernel. This allows for a determination of the dimension of the evaluation of simple functors.\par

As before, we fix a finite group $H$ and we assume that $k$ is algebraically closed (or at least large enough for the group $\Out(H)$). We assume that $V$ is a simple $k\Out(H)$-module, $\sur\pi: k\Out(H)\to E$ is the natural surjection onto $E=\End_k(V)$, and $\tau=\tau_V$ is the trace map on $\End_k(V)$.
In this situation, we have the following.

\result{Theorem} \label{dimension-simple}
Let $V$ be a simple $k\Out(H)$-module, let $\sur\pi: k\Out(H)\to E$ be the natural surjection onto $E=\End_k(V)$, and let $\tau=\tau_V$ be the trace map on $\End_k(V)$. Let ${<}-,-{>}_G$ be the corresponding bilinear form on~$k\sur B(G,H)$, as defined in Section~6. If $G$ is a finite group such that $H\sqsubseteq G$, then
$$\dim(S_{H,V}(G)) = {\rank \,{<}-,-{>}_G\over \dim(V)} \,.$$
\fresult

\pf
As a $k\Out(H)$-module, we have $\End_k(V)\cong m\cdot V$, where $m=\dim(V)$ and $m\cdot V$ denotes the direct sum of $m$ copies of~$V$. Moreover, the functor $W\mapsto S_{H,W}$ is additive in~$W$. Therefore, by Theorem~\ref{B/R}, we obtain
$$k\sur B(G,H)/R(G,H) \cong S_{H,E}(G) \cong m\cdot S_{H,V}(G) \,$$
where $R(G,H)$ is the right kernel of the bilinear form.
Now the dimension of the left hand side is the rank of the bilinear form ${<}-,-{>}_G$ on~$k\sur B(G,H)$. The result follows.
\endpf
 
We now describe a procedure for computing the dimension of $S_{H,V}(G)$, for any finite group $G$. First we need a basis of $k\sur B(G,H)$. We consider all sections $(S,T)$ of $G$ such that $S/T\cong H$, up to $G$-conjugation. Then $\Indinf_{S/T}^G\,\Iso_\sigma$ is a basis element of $k\sur B(G,H)$, where $\sigma:H\to S/T$ is an isomorphism. But we have $\Indinf_{S/T}^G\,\Iso_\sigma=\Indinf_{S/T}^G\,\Iso_\rho$ whenever $\rho=\Conj_g\,\sigma$, where $\Conj_g$ denotes conjugation by some element $g\in N_G(S,T)$.

\result{Lemma} \label{basis}
Let $H$ and $G$ be finite groups. The set of elements of the form $\Indinf_{S/T}^G\,\Iso_\sigma$, where $(S,T)$ runs over all sections of $G$ such that $S/T\cong H$, up to $G$-conjugation, and where $\sigma:H\to S/T$ runs over all isomorphisms up to left composition by $N_G(S,T)$, is a $k$-basis of $k\sur B(G,H)$.
\fresult

\pf
This easily follows from Lemma~\ref{biset} and the analysis above.
\endpf

Now we consider the form ${<}-,-{>}_G$ on $k\sur B(G,H)$. Let $\alpha=\Indinf_{S/T}^G\,\Iso_\sigma$ and $\beta=\Indinf_{J/K}^G\,\Iso_\rho$, where $\sigma:H\to S/T$ and $\rho:H\to J/K$ are isomorphisms. Then
$${<}\alpha,\beta{>}_G=\tau\sur\pi(\alpha\op \beta)$$
and we have $\alpha\op = \Iso_{\sigma^{-1}}\Defres_{S/T}^G$. The generalized Mackey formula (see Proposition~A1 in~\cite{BT1} and Lemma~2.5 in~\cite{BT2})
tells us how to decompose the biset $\Defres_{S/T}^G\,\Indinf_{J/K}^G$ as a sum of transitive bisets (using butterflies, as defined in~\cite{BT2}). Many terms factor through a smaller subquotient and are therefore zero in $k\sur B(G,H)$. We are left with terms involving conjugates of $(J,K)$ which are {\it linked\/} to $(S,T)$ and each of them produces an element of~$\Out(H)$ (see \cite{BT2} for details about linked sections). Computing $\tau\sur\pi$ on such a term $x$ is just computing the character value $\tau_V(x)$. All this is, at least in principle, easy to compute, using some standard computer software of group theory. Note that there is an effective method for computing the evaluation $S_{H,V}(G)$ of a simple functor (see Theorem~4.3.20 in~\cite{Bo3}), but this requires the knowledge of the $k\Out(H)$-module $V$ rather than merely its character.


\section{Semi-simple quotients}
\noindent
We continue with the notation of Section~6 and consider another special case.
We fix a finite group $H$ and
we suppose that $k$ is algebraically closed (or at least large enough for the group $\Out(H)$).
We assume now that the $k$-algebra $E$ is the largest semi-simple quotient of~$k\Out(H)$, that is,
$$E=k\Out(H)/J(k\Out(H))\cong \prod_{i=1}^r \End_k(V_i) \,,$$
where $V_1,\ldots,V_r$ are the simple $k\Out(H)$-modules. The symmetrizing form $\tau:E\to k$ is the sum of all trace maps $\tau=\sum_{i=1}^r \tau_{V_i}$. If $k$ has characteristic~0 (or a prime not dividing $|\Out(H)|$), then $E=k\Out(H)$ and we could also take the ordinary map $\tau$ for group algebras (coefficient of~1). Then, for each finite group $X$, we have a symmetric bilinear form ${<}-,-{>}_X$ on $k\sur B(X,H)$ with kernel $R(X,H)$, and this defines a subfunctor $R(-,H)$ of $k\sur B(-,H)$.

\result{Theorem} \label{semi-simple-quotient}
With the notation above, the following holds.
\begin{enumerate}
\item The semi-simple functor $k\sur B(-,H)/R(-,H)$ is isomorphic to
$$k\sur B(-,H)/R(-,H) \cong \bigoplus_{i=1}^r S_{H,\End_k(V_i)} \cong  \bigoplus_{i=1}^r \;m_i\cdot S_{H,V_i}
\,,$$
where $m_i=\dim(V_i)$ and $m_i\cdot S_{H,V_i}$ denotes the direct sum of $m_i$ copies of~$S_{H,V_i}$.
\item $k\sur B(-,H)/R(-,H)$ is the largest semi-simple quotient of the biset functor $k\sur B(-,H)$. In other words, $R(-,H)$ is the Jacobson radical of $k\sur B(-,H)$.
\end{enumerate}
\fresult

It should be noted that, if each evaluation of a biset functor $F$ is finite-dimensional, the functor $F$ itself may not have finite length. In particular, there are examples of biset functors which have no maximal subfunctors at all. Here, the statement of the theorem shows that $k\sur B(-,H)$ has finitely many maximal subfunctors and that it has a Jacobson radical $R(-,H)$ such that $k\sur B(-,H)/R(-,H)$ has finite length.

\pf
The first part follows immediately from Theorem~\ref{B/R}, together with the following isomorphisms of $k\Out(H)$-modules~:
$$E=k\Out(H)/J(k\Out(H))\cong \bigoplus_{i=1}^r \End_k(V_i)\cong \bigoplus_{i=1}^r \; m_i\cdot V_i \,.$$
To prove part~(2), we assume that a simple functor $S$ is isomorphic to a quotient of~$k\sur B(-,H)$, and we write $S=S_{J,W}$ where $J$ is a group of minimal order such that $S(J)\neq0$ and $W$ is a simple $k\Out(J)$-module. Thus we obtain a non-zero surjective morphism
$$ kB(-,H) \longrightarrow \sur kB(-,H)  \longrightarrow S_{J,W} \,,$$
which must correspond to an element of $S_{J,W}(H)$ by Yoneda's lemma. Therefore $S_{J,W}(H)\neq0$, hence $J\sqsubseteq H$ by Proposition~\ref{parametrization}.
On the other hand, we must have $\sur kB(J,H)\neq0$, because of the surjection
$$\sur kB(J,H)\to S_{J,W}(J)=W\neq0 \,.$$
It follows that $H\sqsubseteq J$ by Lemma~\ref{sqsubset}. Therefore we obtain $J\cong H$.\par

Now $S=S_{H,W}$ is a quotient of $k\sur B(-,H)$. If it was not a quotient of $k\sur B(-,H)/R(-,H)$, we would obtain a semi-simple quotient of $k\sur B(-,H)$ isomorphic to
$$S \oplus k\sur B(-,H)/R(-,H) \;.$$
Then, on evaluation at~$H$, we would obtain a semi-simple quotient of $k\Out(H)$ isomorphic to
$$S_{H,W}(H) \oplus k\Out(H)/R(H,H) = W\oplus k\Out(H)/J(k\Out(H)) \,,$$
because $R(H,H)=\Ker(\sur\pi)=J(k\Out(H))$ by Lemma~\ref{bilinear} and our choice of~$E$. But $k\Out(H)$ cannot have such a semi-simple quotient, by definition of the Jacobson radical. So this is impossible and it follows that $k\sur B(-,H)/R(-,H)$ is the largest semi-simple quotient of $k\sur B(-,H)$.
\endpf

Since the evaluation of a simple functor is a simple module, we obtain the following corollary.

\result{Corollary} \label{semi-simple-evaluation}
With the same notation as above, let $G$ be a finite group.
The $kB(G,G)$-module $k\sur B(G,H)/R(G,H)$ is semi-simple, isomorphic to
$$k\sur B(G,H)/R(G,H) \cong  \bigoplus_{i=1}^r \;m_i\cdot S_{H,V_i}(G) \,,$$
where $m_i=\dim(V_i)$ and $m_i\cdot S_{H,V_i}(G)$ denotes the direct sum of $m_i$ copies of~$S_{H,V_i}(G)$.
\fresult

\pf
This follows immediately from Theorem~\ref{semi-simple-quotient} by taking evaluation at~$G$.
\endpf

We warn the reader that some of the evaluations $S_{H,V_i}(G)$ might be zero. This is an important issue for understanding simple $kB(G,G)$-modules and counting them (see~\cite{BST}).\par

In view of part~(2) of Theorem~\ref{semi-simple-quotient}, the obvious question is to know whether or not $k\sur B(G,H)/R(G,H)$ is the largest semi-simple quotient of~$k\sur B(G,H)$. In other words, when is the inclusion
$$J\big(k\sur B(G,H)\big) \subseteq R(G,H)$$
an equality $J\big(k\sur B(G,H)\big) = R(G,H)$~? We shall see in Section~9 that the answer is negative in general, but positive in some specific cases, for instance when $G$ is abelian.\par

Corollary~\ref{semi-simple-evaluation} produces simple modules for the double Burnside ring $kB(G,G)$ by means of an easily computable quotient $k\sur B(G,H)/R(G,H)$ of a rather straightforward module $k\sur B(G,H)$. This is enough for the determination of the Jacobson radical of $kB(G,G)$, in view of the following result.

\result{Proposition} \label{JBGG}
Let $G$ be a finite group. The Jacobson radical $J\big(kB(G,G)\big)$ of the $k$-algebra $kB(G,G)$
is equal to the kernel of the action of $kB(G,G)$ on the semi-simple module
$$ \bigoplus_{H\sqsubseteq G} k\sur B(G,H)/R(G,H) \,.$$
In other words, every simple $kB(G,G)$-module appears in one of the modules $k\sur B(G,H)/R(G,H)$, for some subquotient $H\sqsubseteq G$.
\fresult

\pf
Let $\alpha\in kB(G,G)$. Then $\alpha\in J\big(kB(G,G)\big)$ if and only if $\alpha$ acts by zero on every simple module $S_{H,V_i}(G)$, where $H\sqsubseteq G$ and $V_1,\ldots,V_r$ are the simple $k\Out(H)$-modules. Equivalently, for every subquotient $H\sqsubseteq G$, $\alpha$ acts by zero on $\bigoplus_{i=1}^r \;m_i\cdot S_{H,V_i}(G)$, that is, on
$k\sur B(G,H)/R(G,H)$ by Corollary~\ref{semi-simple-evaluation}.
\endpf

\result{Remark} \label{multiplicities}
{\rm 
It is interesting to note that $S_{H,V_i}(G)$ (provided it is non-zero) appears at least $m_i=\dim(V_i)$ times as a quotient of $k\sur B(G,H)$, hence also as a quotient of
$$ \bigoplus_{H\sqsubseteq G} k\sur B(G,H)/R(G,H) \,.$$
On the other hand, if $k$ is a splitting field for $kB(G,G)$, the simple module $S_{H,V_i}(G)$ must appear $n_i$ times as a quotient of $kB(G,G)$, where $n_i=\dim(S_{H,V_i}(G))$. The comparison between $m_i$ and $n_i$ is not straightforward. In case $k\sur B(G,H)$ is a generated by a single element (hence a quotient of $kB(G,G)$), then we must have $n_i\geq m_i$. This happens in suitable examples (see Remark~\ref{single} and Proposition~\ref{abelian-etc}). But in general, it seems difficult to have specific information on the dimension of evaluations, apart from the general result of Theorem~\ref{dimension-simple}.
}
\fresult


\section{Further results on evaluations}
\noindent
As in the previous section, $H$ denotes a finite group, $E=k\Out(H)/J\big(k\Out(H)\big)$, and $R(G,H)$ is the kernel of the corresponding bilinear form on~$k\sur B(G,H)$, where $G$ is a finite group (with $H\sqsubseteq G$). In this section, we examine a few cases where the equality $J\big(k\sur B(-,H)\big) = R(-,H)$ remains true on evaluation at~$G$. We first give a sufficient condition.

\result{Proposition} \label{ab=id}
With the notation above, suppose that there exists $\alpha\in kB(H,G)$ and $\beta\in kB(G,H)$ such that $\alpha\beta\equiv \Id\mod kI(H,H)$. Then
$$J\big(k\sur B(G,H)\big) = R(G,H)$$
and the only simple quotients of $k\sur B(G,H)$ (as $kB(G,G)$-module) are the simple modules $S_{H,V}(G)$ (provided they are non-zero), where $V$ is a simple $k\Out(H)$-module.
\fresult
 
\pf
Let $W$ be a simple quotient of $k\sur B(G,H)$ as $kB(G,G)$-module and let $P$ be a minimal group for~$W$. By Proposition~\ref{minimal}, $W=S_{P,V}(G)$ for some simple $k\Out(P)$-module $V$. By Proposition~\ref{quotient}, $H\sqsubseteq P$. Now by assumption
$$k\sur B(G,H)=k\sur B(G,H)\sur\alpha\sur\beta=k\sur B(G,H)\sur\alpha\sur\beta\sur\alpha\sur\beta
\subseteq kB(G,G) \beta\alpha\cdot\sur\beta \subseteq k\sur B(G,H) \,,$$
hence $kB(G,G) \beta\alpha\cdot\sur\beta = k\sur B(G,H)$.
Let $w_0$ be the image of $\sur\beta\in k\sur B(G,H)$ via the quotient map $k\sur B(G,H)\to W$.
Then $\beta\alpha\cdot w_0\neq0$, otherwise the whole image of $k\sur B(G,H)$ would be zero.
But this can be viewed as successive actions of $\alpha$ and $\beta$ in the functor~$S_{P,V}$,
because $W=S_{P,V}(G)$. We deduce that $\alpha\cdot w_0\neq0$ in~$S_{P,V}(H)$.
Therefore $S_{P,V}(H)\neq0$, hence $P\sqsubseteq H$ by Proposition~\ref{parametrization}. It follows that $P\cong H$ and so $W=S_{H,V}(G)$. Therefore, every simple module appearing as a quotient of $k\sur B(G,H)$ must be indexed by~$H$. We are left with the question of the multiplicities.\par

By Corollary~\ref{semi-simple-evaluation}, we know that $k\sur B(G,H)/R(G,H)$ is semi-simple, with $S_{H,V}(G)$ appearing $m$ times, where $m=\dim(V)$. We need to prove that the multiplicity of $S_{H,V}(G)$ as a quotient of $k\sur B(G,H)$ is exactly~$m$. From this, the equality $J\big(k\sur B(G,H)\big) = R(G,H)$ will follow. Now the multiplicity~$\sur m$ of $W=S_{H,V}(G)$ as a quotient of $k\sur B(G,H)$ is equal to
$$\sur m= \dim\Hom_{kB(G,G)}(k\sur B(G,H),W)$$
and there is a $k$-linear map
$$f: V \longrightarrow \Hom_{kB(G,G)}(k\sur B(G,H),W) \,,\quad v\mapsto f(v) \,,$$
where $f(v)$ maps $\sur\gamma\in k\sur B(G,H)$ to $\gamma\cdot v$. This makes sense in the functor $S_{H,V}$, because $S_{H,V}(H)=V$ and $S_{H,V}(G)=W$, and it only depends on~$\sur\gamma$, because $S_{H,V}$ vanishes on proper subquotients of~$H$.
The map $f$ is non-zero because there is some $\gamma\in kB(G,H)$ and $v\in V$ such that $\gamma\cdot v\neq0$, by Proposition~\ref{generated}. We shall prove that $f$ is an isomorphism, so that $m=\dim(V)$ must be equal to~$\sur m$, as required.

Now the map $f$ is $k\Out(H)$-linear, with respect to the left action on the target given by
$$(\sur\phi\cdot q)(\sur\gamma)=q(\sur\gamma\sur\phi) \,,$$
where $\sur\phi\in k\Out(H)$, $q\in\Hom_{kB(G,G)}(k\sur B(G,H),W)$, and $\sur\gamma\in k\sur B(G,H)$.
By simplicity of $V$ as a $k\Out(H)$-module, the map $f$ is injective. To prove the surjectivity of~$f$, we let $q\in\Hom_{kB(G,G)}(k\sur B(G,H),W)$. By our assumption, we have
$$q(\sur\gamma)=q(\sur\gamma\,\sur\alpha\sur\beta)=\gamma\alpha\, q(\sur\beta)$$
by $kB(G,G)$-linearity of~$q$. Therefore $q(\sur\gamma)=\gamma\cdot v$, where $v=\alpha\,q(\sur\beta)$, and so $q=f(v)$, proving the surjectivity.
\endpf

\result{Remark} \label{char}
{\rm If $k$ has characteristic zero or prime to $|\Out(H)|$, then we can avoid the argument of the second part of the proof and prove directly that $\sur m=m$, as follows. We already know that $k\sur B(G,H)$ has a semi-simple quotient with $S_{H,V}(G)$ appearing $m$ times, where $m=\dim(V)$. If the multiplicity of $S_{H,V}(G)$ is $\sur m\geq m$, then $S_{H,V}$ appears $\sur m$ times as a composition factor of~$k\sur B(-,H)$, by repeated applications of Proposition~\ref{subquotients}. Then, on evaluation at~$H$, the module $V=S_{H,V}(H)$ appears at least $\sur m$ times as a composition factor of~$k\sur B(H,H)=k\Out(H)$. But since $k\Out(H)$ is semi-simple by Maschke's theorem, the module $V$ appears exactly $m=\dim(V)$ times as a composition factor of~$k\Out(H)$. Therefore $\sur m\leq m$. 
}
\fresult

\result{Remark} \label{single}
{\rm 
The assumption of Proposition~\ref{ab=id} implies that $k\sur B(G,H)$ is generated by a single element as a $kB(G,G)$-module, because $kB(G,G) \beta\alpha\cdot\sur\beta = k\sur B(G,H)$, as we have seen in the proof. However, there are many examples where $k\sur B(G,H)$ is not a cyclic $kB(G,G)$-module (e.g. Example~\ref{example-A5} below). It seems to be an interesting question to understand when this happens. It clearly depends on the various ways $H$ is realized as a subquotient of the group $G$.
}
\fresult

We now show that the situation of Proposition~\ref{ab=id} is rather common, so that the equality $J\big(k\sur B(G,H)\big) = R(G,H)$ often occurs.

\result{Proposition} \label{abelian-etc}
The assumption of Proposition~\ref{ab=id} is satisfied in each of the following cases.
\begin{enumerate}
\item $H$ is isomorphic to a quotient of $G$.
\item $G$ is abelian.
\item $H$ is isomorphic to a subgroup $Z$ of $G$ such that $N_G(Z)=ZC_G(Z)$, and $|N_G(Z):Z|$ is non-zero in~$k$.
\item $H$ is isomorphic to a central subgroup $Z$ of $G$, and $|G:Z|$ is non-zero in~$k$.

\end{enumerate}
\fresult

\pf
(1) If $H=G/N$, then $\Def_{G/N}^G\,\Inf_{G/N}^G =\Id_{G/N}$.\par

(2) Any subgroup of an abelian group $G$ is isomorphic to a quotient of~$G$, by using the isomorphism between $G$ and its dual. It follows that any subquotient of $G$ is isomorphic to a quotient of~$G$. Thus part~(1) applies.\par

(4) is a special case of (3). Note that, in characteristic zero, (2) can also be proved by using~(4).\par

We are left with a proof of~(3). Let $\alpha=|N_G(Z):Z|^{-1} \Res_Z^G$ and $\beta=\Ind_Z^G$. By the Mackey formula, we have
\begin{eqnarray*}
\Res_Z^G\,\Ind_Z^G
&=&\sum_{g\in[Z\backslash G/Z]} \Ind_{Z\cap \ls gZ}^Z\,\Res_{Z\cap \ls gZ}^{\ls gZ} \,\Conj_g \\
&\equiv&\sum_{g\in[N_G(Z)/Z]} \Ind_Z^Z\,\Res_Z^Z\,\Conj_g \mod I(Z,Z) \\
&\equiv& |N_G(Z):Z| \cdot \Id  \mod I(Z,Z) \,,
\end{eqnarray*}
because $\Conj_g$ is the identity for $g\in ZC_G(Z)=N_G(Z)$.
It follows that $\alpha\beta\equiv \Id\mod kI(Z,Z)$.
\endpf

\result{Remark} \label{more-general}
{\rm 
Statement (3) is in fact a special case of a more general, but more technical, result. Suppose that $H$ is isomorphic to a subquotient $S/T$ of~$G$, where $(S,T)$ is a section of $G$ such that $T\leq \Phi(S)$ and the image of $N_G(S,T)$ in $\Out(S)$ is trivial. Suppose also that $|N_G(S,T):S|$ is non-zero in~$k$. Then the assumption of Proposition~\ref{ab=id} is satisfied. The proof uses the idempotents $e_S^G$ of the ordinary Burnside ring $B(G)$ and their images $\widetilde{e_S^G}\in B(G,G)$ defined in Section~2.5 of~\cite{Bo3}.
We let
$$\alpha = |N_G(S,T):S|^{-1} \Defres_{S/T}^G\, \widetilde{e_S^G} \quad \text{and} \quad
\beta=\widetilde{e_S^G}\, \Indinf_{S/T}^G \,.$$
Then the actual computation is similar in spirit to the one used in the proof above, but more involved. One needs the generalized Mackey formula (Proposition~A1 in~\cite{BT1} and Lemma~2.5 in~\cite{BT2}), which tells us how to decompose the biset $\Defres_{S/T}^G\,\Indinf_{S/T}^G$ as a sum of transitive bisets, using butterflies, as defined in~\cite{BT2}. The condition that $T$ is contained in the Frattini subgroup $\Phi(S)$ implies that $\widetilde{e_S^G}\, \Indinf_{S/T}^G=\Indinf_{S/T}^G \, \widetilde{e_{S/T}^{S/T}}$ and this is used to show that many terms in the the sum lie in fact in~$I(S/T,S/T)$. The only remaining terms are conjugations by elements $g\in N_G(S,T)$. The assumption on the action of this group ensures that such a conjugation is the identity and the result follows.
}
\fresult

\result{Remark} \label{expansive}
{\rm 
Yet another case where we obtain the equality $J\big(k\sur B(G,H)\big) = R(G,H)$ appears when $H=N_G(T)/T$ and $T$ is an {\it expansive\/} subgroup of~$G$, as defined in Section~6.4 of~\cite{Bo3}. We consider the elements
$$\alpha = f_1^{N_G(T)/T}\,\Defres_{N_G(T)/T}^G \quad \text{and} \quad
\beta= \Indinf_{S/T}^G \, f_1^{N_G(T)/T} \,,$$
where $f_1^{N_G(T)/T}$ is also defined in Section~6.4 of~\cite{Bo3}.
Then $\alpha$ and $\beta$ satisfy the property $\alpha\beta\equiv \Id\mod kI(H,H)$ and Proposition~\ref{ab=id} applies again.
}
\fresult

\bigskip
We have seen various cases where $J\big(k\sur B(G,H)\big) = R(G,H)$, but we shall see in Section~13 several examples where this is not so.


\section{The case of the trivial group}
\noindent
We assume now that the finite group $H$ is trivial. Then $k\sur B(G,1)=kB(G,1)=kB(G)$ for any finite group $G$, where $B(G)$ denotes the ordinary Burnside ring of~$G$. Viewed as biset functors, we have
$$k\sur B(-,1)=kB(-,1)=kB(-)=L_{1,k}=\sur L_{1,k} \,.$$
Moreover, this has a unique simple quotient $S_{1,k}$ by Lemma~\ref{JGW} (or also by Theorem~\ref{semi-simple-quotient}). For simplicity, we assume that $k$ is a field of characteristic zero. In this special case, our bilinear form on $k\sur B(G,1)=kB(G)$ was already considered in Section~7.2 of~\cite{Bo1}. The kernel $R(G,1)$ of the bilinear form is equal to the kernel of the surjective map
$$q: kB(G)\to kR_\Q(G) \,,$$
where $R_\Q$ is the functor of rational representations. Here $q$ denotes the $k$-linear extension of the natural homomorphism $q:B(G)\to R_\Q(G)$ mapping a $G$-set~$X$ to the permutation $\Q G$-module $q(X)=\Q X$. It follows that $S_{1,k}\cong kR_\Q$. In particular the dimension of $S_{1,k}(G)$ is the number of conjugacy classes of cyclic subgroups of~$G$.\par

Now we extend further the analysis by introducing the two-sided ideal $I(G)$ of $kB(G,G)$ generated by all the $(G,G)$-bisets which factor through the trivial group. 
The method finds its origin in the proof of Proposition~6.1.5 in~\cite{Bo3}. Thus $I(G)$ is generated as a $k$-vector space by the $(G,G)$-bisets
$$(G\times G)\big/(A\times B)=\Indinf_{A/A}^G \, \Iso_\sigma \, \Defres_{B/B}^G \,,$$
where $A$ and $B$ are subgroups of $G$ and $\sigma: B/B \to A/A$ is the obvious unique isomorphism (which we shall ignore in the sequel, for simplicity). It is useful to view $G/A$ as a $(G,1)$-biset and $B\backslash G$ as a $(1,G)$-biset, that is,
$$G/A=\Indinf_{A/A}^G \quad\text{and}\quad B\backslash G =(G/B)\op=\Defres_{B/B}^G \,.$$
With this point of view, we obtain
$$(G\times G)\big/(A\times B)= (G/A) \cdot (G/B)\op \,.$$
More generally, any $\alpha\in kB(G)$ can be viewed as an element of $kB(G,1)$ and then
$$\alpha\beta\op\in I(G) \,, \quad\text{for any } \;\alpha,\beta\in kB(G) \,.$$
This applies in particular to the primitive idempotents $e_A^G$ of the Burnside ring $kB(G)$, which form a $k$-basis of $kB(G)$, and we obtain generators $e_A^G \,(e_B^G)\op$ of~$I(G)$ as a vector space. Here $A$ and $B$ run over all subgroups of~$G$ up to conjugation.\par

The action of $I(G)$ on $S_{1,k}(G)=kR_\Q(G)$ is described in the following result.

\result{Lemma} \label{action}
Let $k$ be a field of characteristic zero and let $I(G)$ be the two-sided ideal of $kB(G,G)$ defined above.
\begin{enumerate}
\item If $H$ is a non-trivial finite group, the action of $I(G)$ on~$S_{H,V}(G)$ is zero, for any $k\Out(H)$-module~$V$.
\item Let $\alpha,\beta\in kB(G)$ and let $M$ be a $\Q G$-module. The action of $\alpha\beta\op$ on~$[M]\in kR_\Q(G)$ is given by~:
$$\alpha\beta\op \cdot [M] =  \big(q(\beta)\mid M\big)_G \cdot q(\alpha) \,,$$
where $\big(-\mid-\big)_G$ denotes the ordinary scalar product of $\Q G$-modules.
\item If either $A$ or $B$ is non-cyclic, then the action of $e_A^G \,(e_B^G)\op$ on~$kR_\Q(G)$ is zero.
\item If $A$ and $B$ run over the set of cyclic subgroups up to $G$-conjugation, then the images of $e_A^G \,(e_B^G)\op$ in~$\End_k\big(kR_\Q(G)\big)$ are $k$-linearly independent. So are the images of $G/A \cdot(G/B)\op$, where $A$ and $B$ run over the same set.
\end{enumerate}
\fresult

\pf
(1) The action of an element of $I(G)$ factors through the trivial group.
But $S_{H,V}(1)=0$ if $H\neq1$ since $H$ has minimal order such that $S_{H,V}(H)\neq0$.\par

(2) Taking basis elements, we can assume that $\alpha=G/A=\Indinf_{A/A}^G$ and $\beta=G/B$, so $\beta\op=\Defres_{B/B}^G$. Then
$$G/A\cdot(G/B)\op(M) =\Indinf_{A/A}^G \, \Defres_{B/B}^G(M) =\Indinf_{A/A}^G(M^B)$$
because the deflation of a $\Q B$-module is obtained by taking $B$-fixed points (actually cofixed points, but this is the same in characteristic~0). Then $M^B$ is just the direct sum of $\dim(M^B)$ copies of $\Q$ (a module for the trivial group), and then this is induced from $A$ to $G$. So we obtain
\begin{eqnarray*}
G/A\cdot(G/B)\op(M) &=&\dim(M^B)\cdot \Q[G/A]\\
&=&\big(\Q\mid \Res_B^G(M)\big)_B\cdot q(G/A)\\
&=&\big(\Q[G/B]\mid M\big)_G\cdot q(G/A)\\
&=&\big(q(G/B)\mid M\big)_G\cdot q(G/A) \,,
\end{eqnarray*}
as was to be shown.\par

(3) If $A$ is non-cyclic, then it is well-known that $q(e_A^G)=0$ (because the restriction of $e_A^G$ to any cyclic subgroup of $G$ is zero, hence the same holds for $q(e_A^G)$). So the action of $e_A^G \,(e_B^G)\op$ is zero by~(2). The same holds if $B$ is non-cyclic, because $q(e_B^G)=0$ in that case.\par

(4) Suppose that $\sum_{A,B}\lambda_{A,B}\,G/A \cdot(G/B)\op$ acts by zero on $kR_\Q(G)$, where $A$ and $B$ run over the set of cyclic subgroups up to $G$-conjugation and where $\lambda_{A,B}\in k$. Then by (2), for every $x\in kR_\Q(G)$, we have
$$\sum_A\sum_B\lambda_{A,B}\,\big(\Q[G/B]\mid x\big) \cdot \Q[G/A] =0 \,.$$
Since the modules $\Q[G/A]$ form a $k$-basis of $kR_\Q(G)$ by Artin's induction theorem, we obtain
$$\big(\sum_B\lambda_{A,B}\,\Q[G/B]\mid x\big)=0 \,.$$
But this holds for all $x$ and the scalar product is non-degenerate. Therefore
$$\sum_B\lambda_{A,B}\,\Q[G/B]=0\,.$$
Again the modules $\Q[G/B]$ form a basis, so $\lambda_{A,B}=0$. This shows one of the statement. The other statement follows because the idempotents $e_A^G$ generate the same subspace of the Burnside ring $kB(G)$ as the elements $G/A$ (where $A$ runs over cyclic subgroups up to conjugation).
\endpf

Now we can prove that a large part of the ideal $I(G)$ lies in the Jacobson radical $J(kB(G,G))$.

\result{Theorem} \label{trivial-group}
Let $k$ be a field of characteristic zero and let $I(G)$ be the ideal of $kB(G,G)$ defined above. Let $I_c(G)$ be the $k$-subspace of $I(G)$ generated by all elements $(G/A)\cdot(G/B)\op$ such that $A$ and $B$ are cyclic subgroups of~$G$.
\begin{enumerate}
\item $kB(G,G)=I_c(G)\oplus M$, where $M$ denotes the maximal two-sided ideal of $kB(G,G)$ which is the kernel of  the action of $kB(G,G)$ on $S_{1,k}(G)=kR_\Q(G)$.
\item $I(G)=I_c(G)\oplus \big(I(G)\cap J(kB(G,G))\big)$.
\item $\dim\big(I(G)\cap J(kB(G,G))\big) = b(G)^2-c(G)^2$, where $b(G)$ is the number of conjugacy classes of subgroups of $G$ and $c(G)$ is the number of conjugacy classes of cyclic subgroups of~$G$.
\end{enumerate}
\fresult

\pf
(1) Since $kR_\Q(G)=S_{1,k}(G)$ is a simple $kB(G,G)$-module, the natural $k$-algebra homomorphism
$$r:kB(G,G) \longrightarrow \End_k(kR_\Q(G))$$
is surjective by the density theorem. Part~(4) of Lemma~\ref{action} shows that there are $c(G)^2$ elements of $I_c(G)$ whose images under $r$ are $k$-linearly independent. But $kR_\Q(G)$ has dimension~$c(G)$ (by Artin's induction theorem), so $\dim\big(\End_k(kR_\Q(G))\big)=c(G)^2$. It follows that the restriction of $r$ to $I_c(G)$ is an isomorphism of $k$-vector spaces and (1) follows (because $\Ker(r)=M$).\par

(2) Let $I'(G)$ be the $k$-subspace of $I(G)$ generated by all elements $e_A^G \,(e_B^G)\op$, where either $A$ or $B$ is non-cyclic. Then $I(G)=I_c(G)\oplus I'(G)$ because the elements $e_A^G \,(e_B^G)\op$ where  $A$ and $B$ are both cyclic generate $I_c(G)$. By parts (1) and (3) of Lemma~\ref{action}, $I'(G)$ acts by zero on every simple $kB(G,G)$-module, hence $I'(G)\subseteq J(kB(G,G))$. Since $I_c(G)$ acts faithfully on one simple module (by part~(4) of Lemma~\ref{action} ), it follows that $I'(G)=I(G)\cap J(kB(G,G))$.\par

(3) It is clear that $\dim(I(G))=b(G)^2$ and $\dim(I_c(G))=c(G)^2$.
\endpf

We now recover a result proved in Section~6.1 of~\cite{Bo3} (but the result in~\cite{Bo3} is more precise).

\result{Corollary} \label{not-semi-simple} 
If $k$ is a field of characteristic zero and if $G$ is a non cyclic group, then $kB(G,G)$ is not semi-simple.
\fresult

\pf
With the notation of Theorem~\ref{trivial-group}, we have $b(G)>c(G)$ because $G$ is not cyclic.
Therefore $I(G)\cap J(kB(G,G))\neq0$ and $kB(G,G)$ is not semi-simple.
\endpf


\section{Left ideals in the double Burnside ring}
\noindent
In order to understand the ring structure of the double Burnside ring $kB(G,G)$, it is useful to have information on some naturally defined left ideals associated to sections of~$G$. Many of the previous results used a finite group $H$ and the space $kB(G,H)$, or its quotient $k\sur B(G,H)$. Now we work with an incarnation of~$H$ as a subquotient of~$G$ by fixing a section $(P,Q)$ of~$G$ such that $P/Q= H$.
Then the group $N_G(P,Q)$ comes into play, or more precisely the image $\Gamma_G(P,Q)$ of $N_G(P,Q)$ in $\Out(P/Q)$.\par

We compose with $ \Defres_{P/Q}^G$ in order to obtain elements of $kB(G,G)$, that is, we define
\begin{eqnarray*}
K_{(P,Q)}&=&  \{\, \alpha\cdot \Defres_{P/Q}^G \,|\, \alpha\in kB(G,P/Q) \,\} \\
K_{<(P,Q)}&=& \{\, \alpha \cdot \Defres_{P/Q}^G \,|\, \alpha\in kI(G,P/Q) \,\} \\
\sur K_{(P,Q)}&=&  K_{(P,Q)}/K_{<(P,Q)} \,.
\end{eqnarray*}
Clearly,  $K_{(P,Q)}$ and $K_{<(P,Q)}$ are left ideals of $kB(G,G)$ and hence $\sur K_{(P,Q)}$ is a $kB(G,G)$-module. Moreover, $K_{(P,Q)}$ and $\sur K_{(P,Q)}$ only depend on the $G$-conjugacy class of the section~$(P,Q)$. Using Lemma~\ref{basis}, one can find a $k$-basis of $\sur K_{(P,Q)}$ consisting of elements of the form $\Indinf_{J/K}^G\,\Iso_\sigma\,\Defres_{P/Q}^G$, where $(J,K)$ runs over sections of $G$ and $\sigma:P/Q\to J/K$ is an isomorphism, up to left composition by $N_G(J,K)$ and right composition by $N_G(P,Q)$.\par

It is possible to filter the whole ring $kB(G,G)$ by left ideals in such a way that each successive quotient is isomorphic to $\sur K_{(P,Q)}$ for some section $(P,Q)$, as follows. We consider the set $\cal X$ of all conjugacy classes of sections of~$G$, and we let
$$\emptyset={\cal X}_0\subset {\cal X}_1 \subset \ldots \subset {\cal X}_{N-1}\subset {\cal X}_N={\cal X}$$
be a filtration of $\cal X$ such that, for each~$1\leq i\leq N$,
$${\cal X}_i={\cal X}_{i-1}\cup \{[(P_i,Q_i)]\} \,,$$
where $[(P_i,Q_i)]$ is the conjugacy class of a section~$(P_i,Q_i)$ such that $|X/Y|\leq |P_i/Q_i|$ for all $[(X,Y)]\in{\cal X}_{i-1}$. It is clear that such a filtration always exists. Then we define
$$K_{{\cal X}_i}=\sum_{[(X,Y)]\in{\cal X}_i} K_{(X,Y)} \,,$$
and we obtain a filtration
$$0=K_{{\cal X}_0}\subset K_{{\cal X}_1} \subset \ldots \subset K_{{\cal X}_{N-1}}
\subset K_{{\cal X}_N}=kB(G,G) \,.$$
The following result is easy and its proof is left to the reader.

\result{Lemma} \label{filtration} 
With the notation above, $K_{{\cal X}_i}/K_{{\cal X}_{i-1}} \cong \sur K_{(P_i,Q_i)}$.
\fresult

This result shows that the module structure of $\sur K_{(P,Q)}$ is relevant for the understanding of the ring structure of $kB(G,G)$. We shall describe a semi-simple quotient of $\sur K_{(P,Q)}$ when $k$ has characteristic zero (but this simplifying assumption could be dropped with a little more work). We first need the following description of $\sur K_{(P,Q)}$.

\result{Proposition} \label{left-ideal}
Let $G$ be a finite group and let $(P,Q)$ be a section of~$G$.
\begin{enumerate}
\item The homomorphism $kB(G,P/Q)\to K_{(P,Q)}$, $\alpha\mapsto\alpha\cdot \Defres_{P/Q}^G$, induces an isomorphism
$$\mu:kB(G,P/Q)_{N_G(P,Q)}\flh{\sim}{} K_{(P,Q)}$$
where the left hand side denotes the quotient of $kB(G,P/Q)$ by the right action of $N_G(P,Q)$ by conjugation (cofixed points).
\item The isomorphism $\mu$ induces an isomorphism
$$\sur\mu:k\sur B(G,P/Q)_{\Gamma_G(P,Q)}=k\sur B(G,P/Q)_{N_G(P,Q)}\flh{\sim}{} \sur K_{(P,Q)}$$
where $\Gamma_G(P,Q)$ is the image of $N_G(P,Q)$ in $\Out(P/Q)$.
\item The isomorphism $\sur\mu$ induces an isomorphism
$$k\sur B(G,P/Q)\otimes_{k\Out(P/Q)} k[\Out(P/Q)/\Gamma_G(P,Q)] \flh{\sim}{} \sur K_{(P,Q)} \,.$$
\end{enumerate}
\fresult

\pf
(1) If $g\in N_G(P,Q)$, then
$$\alpha\,\Conj_g\, \Defres_{P/Q}^G=\alpha\, \Defres_{P/Q}^G\,\Conj_g=\alpha\, \Defres_{P/Q}^G\,,$$
so that the surjective homomorphism $kB(G,P/Q)\to K_{(P,Q)}$ passes to the quotient $kB(G,P/Q)_{N_G(P,Q)}$ of~$kB(G,P/Q)$.
Then the induced map $\mu$ is easily seen to be injective (by Lemma~\ref{biset}).\par

(2) Taking cofixed points is right exact, so we have a commutative diagram
$$\begin{array}{ccccc}
kB(G,P/Q) & \longrightarrow & k\sur B(G,P/Q) & \longrightarrow & 0 \\
\left\downarrow\vbox to 6mm{}\right.\rlap{$\displaystyle{}$} &&
\left\downarrow\vbox to 6mm{}\right.\rlap{$\displaystyle{}$} && \\
kB(G,P/Q)_{N_G(P,Q)} & \longrightarrow & k\sur B(G,P/Q)_{N_G(P,Q)}& \longrightarrow & 0 \\ \\
\left\downarrow\vbox to 6mm{}\right.\rlap{$\displaystyle{\cong}$} &&&& \\
K_{(P,Q)} & \longrightarrow & \sur K_{(P,Q)} & \longrightarrow & 0 \\
\end{array}
$$
The kernel of the map in the first row is $I(G,P/Q)$. The kernel of the map in the second row is the image $I'$ of $I(G,P/Q)$ in $kB(G,P/Q)_{N_G(P,Q)}$. Now the image of $I'$ in $K_{(P,Q)}$ is the image of $I(G,P/Q)$ in $K_{(P,Q)}$, namely $K_{<(P,Q)}$ by definition. Therefore the vertical isomorphism induces an isomorphism
$$kB(G,P/Q)_{N_G(P,Q)} \big/ I' \, \flh{\sim}{} \, K_{(P,Q)} \big/ K_{<(P,Q)} \,,$$
that is, an isomorphism 
$$k\sur B(G,P/Q)_{N_G(P,Q)} \, \flh{\sim}{} \, \sur K_{(P,Q)} \,,$$
as required. The equality  $k\sur B(G,P/Q)_{N_G(P,Q)}=k\sur B(G,P/Q)_{\Gamma_G(P,Q)}$ is clear, since $N_G(P,Q)$ acts via its image $\Gamma_G(P,Q)$.\par

(3) For the right $k\Out(P/Q)$-module $k\sur B(G,P/Q)$, taking cofixed points under $\Gamma_G(P,Q)$ is the same as tensoring with $k[\Out(P/Q)/\Gamma_G(P,Q)]$.
\endpf

\result{Remark} \label{over-Z} 
{\rm Actually, we do not need to extend scalars to~$k$ in Proposition~\ref{left-ideal}. The result holds for similarly defined ideals in $B(G,G)$.
}
\fresult

For simplicity, we assume now that the field $k$ has characteristic zero and is algebraically closed (or at least large enough for the group $\Out(P/Q)$). Then fixed points and cofixed points are isomorphic (by Maschke's theorem), so we consider the subspace of fixed points $k\sur B(G,P/Q)^{\Gamma_G(P,Q)}$. As in Section~8, we consider the bilinear form ${<}-,-{>}_G$ on $k\sur B(G,P/Q)$ associated with the usual symmetrizing form $\tau$ on the group algebra $k\Out(P/Q)$. The bilinear form ${<}-,-{>}_G$ restricts to the subspace $k\sur B(G,P/Q)^{\Gamma_G(P,Q)}$, hence defines a bilinear form, still denoted ${<}-,-{>}_G$, on $\sur K_{(P,Q)}$ (by part~(2) of Proposition~\ref{left-ideal}). We let $R'$ be the kernel of this bilinear form on $\sur K_{(P,Q)}$.

Now we can prove a result analogous to Corollary~\ref{semi-simple-evaluation}, where we replace $k\sur B(G,H)$ by its `incarnation' $\sur K_{(P,Q)}$ inside the double Burnside ring. It turns out that the semi-simple module $k\Out(H)$ has to be replaced by the module $W=k[\Out(P/Q)/\Gamma_G(P,Q)]$.

\result{Theorem} \label{K-modulo-kernel}
Assume $k$ is algebraically closed of characteristic zero.
Let $G$ be a finite group, let $(P,Q)$ be a section of~$G$, let $K_{(P,Q)}$ and $K_{<(P,Q)}$ be the ideals of $kB(G,G)$ defined above, and let $\sur K_{(P,Q)}=K_{(P,Q)}/K_{<(P,Q)}$.
Let also $R'$ be the kernel of the bilinear form ${<}-,-{>}_G$ on $\sur K_{(P,Q)}$.
Then $\sur K_{(P,Q)}\big/R'$ is a semi-simple $kB(G,G)$-module, isomorphic to
$$\sur K_{(P,Q)} \big/R' \cong S_{P/Q,W}(G) \cong  \bigoplus_{j=1}^q \;m_j\cdot S_{P/Q,W_j}(G) \,,$$
where $W=k[\Out(P/Q)/\Gamma_G(P,Q)]$ and $W=\bigoplus_{j=1}^q \;m_j\cdot W_j$ is the decomposition of the $k\Out(P/Q)$-module $W$ as a direct sum of simple modules.
\fresult

\pf
We write $H=P/Q$. Since $k$ has characteristic zero, $J(k\Out(H))=0$. By Theorem~\ref{semi-simple-quotient}, we have a short exact sequence of biset functors
$$0\longrightarrow R(-,H) \longrightarrow k\sur B(-,H) \longrightarrow S_{H,k\Out(H)} \longrightarrow 0
\,,$$
with a semi-simple right-hand side
$$S_{H,k\Out(H)}\cong\bigoplus_{i=1}^r \;n_i\cdot S_{H,V_i} \,,$$
where $k\Out(H)=\bigoplus_{i=1}^r \;n_i\cdot V_i$ is the decomposition of $k\Out(H)$ as a direct sum of simple modules.
Since $k$ has characteristic zero, tensoring with $W$ is exact and we obtain a short exact sequence of biset functors
\begin{eqnarray*}
0\longrightarrow R(-,H)\otimes_{k\Out(H)}W &\longrightarrow& k\sur B(-,H)\otimes_{k\Out(H)}W \\
&\qquad\longrightarrow & \quad S_{H,k\Out(H)}\otimes_{k\Out(H)}W \longrightarrow 0 \,.
\end{eqnarray*}
Evaluating at $G$, the middle term is isomorphic to $\sur K_{(P,Q)}$, by Proposition~\ref{left-ideal}. We claim that the kernel on the left hand side (evaluated at~$G$) is isomorphic to the kernel $R'$ of the bilinear form. Then it follows that the right hand side (evaluated at~$G$) is
$$S_{H,k\Out(H)}(G)\otimes_{k\Out(H)}W \cong \sur K_{(P,Q)} \big/R' \,.$$

We shall return to this at the end of the proof, but we first prove the claim. We note that tensoring with $W$ is the same as taking (co)fixed points under~$\Gamma$, where $\Gamma=\Gamma_G(P,Q)$. Moreover
$$k\sur B(G,H)=k\sur B(G,H)^\Gamma \oplus U \,,$$
where $U$ is the kernel of (right) multiplication by $e_\Gamma={1\over|\Gamma|}\sum_{\gamma\in\Gamma}\gamma$. This direct sum is orthogonal with respect to the bilinear form, because, if $\sur\phi\in k\sur B(G,H)^\Gamma$ and $\sur\alpha\in U$, we have
$${<}\sur\phi,\sur\alpha{>}_G={<}\sur\phi\,e_\Gamma,\sur\alpha{>}_G
={<}\sur\phi,\sur\alpha\,e_\Gamma{>}_G=0 \,,$$
using part~(5) of Lemma~\ref{bilinear} and the fact that $e_\Gamma\op=e_\Gamma$. Since the direct sum is orthogonal, the kernel of the bilinear form decomposes as
$$R(G,H)=(R(G,H)\cap k\sur B(G,H)^\Gamma) \oplus (R(G,H)\cap U) 
=R(G,H)^\Gamma \oplus (R(G,H)\cap U) \,.$$
Thus the kernel of the bilinear form restricted to $k\sur B(G,H)^\Gamma$ is $R(G,H)^\Gamma$. In other words, in terms of tensor products, the kernel of the bilinear form restricted to $k\sur B(G,H)\otimes_{k\Out(H)}W$ is $R(G,H)\otimes_{k\Out(H)}W$. This completes the proof of the claim.\par

We now return to the isomorphism
$$S_{H,k\Out(H)}(G)\otimes_{k\Out(H)}W \cong \sur K_{(P,Q)} \big/R'$$
and we analyze the left-hand side. We do this by using again the biset functor $S_{H,k\Out(H)}\otimes_{k\Out(H)}W$. Since tensoring with $W$ is taking cofixed points, hence a quotient, $S_{H,k\Out(H)}\otimes_{k\Out(H)}W$ is a quotient of $S_{H,k\Out(H)}$ and is therefore semi-simple again. So we only need to know the simple functors which occur and their multiplicity. This can be achieved by evaluating at~$H$ because all the simple functors which appear are indexed by~$H$. So we evaluate at $H$ and we obtain
$$S_{H,k\Out(H)}(H)\otimes_{k\Out(H)}W =k\Out(H)\otimes_{k\Out(H)}W\cong W =S_{H,W}(H) \,.$$
This forces the isomorphism
$$S_{H,k\Out(H)}\otimes_{k\Out(H)}W\cong S_{H,W}=\bigoplus_{j=1}^q \;m_j\cdot S_{H,W_j}$$
because we know that the left-hand side must be a direct sum of simple functors indexed by~$H$.
Evaluating now at $G$, we obtain the desired isomorphism
$$ \sur K_{(P,Q)} \big/R'\cong S_{H,k\Out(H)}(G)\otimes_{k\Out(H)}W \cong S_{H,W}(G) \,,$$
and the proof is complete.
\endpf

\result{Corollary} \label{abelian}
With the same assumptions, suppose moreover that $G$ is abelian.
Then $R'(G,P/Q)=J(\sur K_{(P,Q)})$ and the corresponding semi-simple quotient is
$$\sur K_{(P,Q)}\big/R'(G,P/Q)\cong k\sur B(G,P/Q)\big/R(G,P/Q)
\cong \bigoplus_{i=1}^r \;m_i\cdot S_{P/Q,V_i}(G) \,,$$
where $k\Out(P/Q) =\bigoplus_{i=1}^r \;m_i\cdot V_i$ is the decomposition of $k\Out(P/Q)$ as a direct sum of simple modules (so $m_i=\dim(V_i)$).
\fresult

\pf
Since $G$ is abelian, the action of $N_G(P,Q)$ on $P/Q$ is trivial, that is, $\Gamma_G(P,Q)=1$. Therefore
$$\sur K_{(P,Q)}\cong k\sur B(G,P/Q) \quad \text{and} \quad
\sur K_{(P,Q)}\big/R'(G,P/Q)\cong k\sur B(G,P/Q)\big/R(G,P/Q) \,.$$
Now by Propositions~\ref{ab=id} and~\ref{abelian-etc} (using again the fact that $G$ is abelian),
$R(G,P/Q)$ is the Jacobson radical of $k\sur B(G,P/Q)$, and hence $R'(G,P/Q)$ is the
Jacobson radical of $\sur K_{(P,Q)}$. The explicit description of the semi-simple quotient follows now from Theorem~\ref{K-modulo-kernel}.
\endpf


\section{The bifree double Burnside ring}
\noindent
In this section, we explain briefly that the results of this paper can be adapted to the case of the bifree double Burnside ring $kA(G,G)$. We then recover some special cases of results of Webb~\cite{We2}. Let $kA(G,H)$ be the subspace of $kB(G,H)$ generated by all the $(G,H)$-bisets which are free on the left and on the right. In other words, we do not allow for inflation and deflation maps, so that $kA(G,H)$ is generated by all bisets of the form
$\Ind_A^G\, \Iso_{\sigma}\,\Res_S^H$, where $A$ is a subgroup of $G$, $S$ is a subgroup of $H$, and $\sigma:S\to A$ is a group isomorphism. In particular $kA(G,G)$ is a subring of $kB(G,G)$.\par

The version used in homotopy theory is the subring of $kB(G,G)$ generated by all the $(G,G)$-bisets which are free on one side, so that deflation is not allowed but inflation is still defined. Unfortunately, our methods do not seem to apply to this subring, because we use symmetry, in particular in the definition of the bilinear forms where opposite bisets are used (and the opposite of inflation is deflation). Therefore, we only consider here the bifree version $kA(G,G)$. We could work with a more general version of the double Burnside ring by using, as in~\cite{We2}, a family of groups $\cal X$ (respectively~$\cal Y$) allowed as normal subgroups involved in deflation (respectively inflation), and our methods would work in case $\cal X = \cal Y$. But for simplicity, we only consider the case ${\cal X} = {\cal Y} = \{1\}$, corresponding to bifree bisets.\par

By restricting to bifree bisets only, the category of biset functors becomes the category of global Mackey functors. All the methods of this paper can be easily adapted and therefore our results hold in this different context. Throughout our arguments, it suffices to replace the sections of a group $G$ by the subgroups of~$G$, so that the relation $H\sqsubseteq G$ now means that $H$ is isomorphic to a subgroup of~$G$. The definition of the standard quotient $k\sur A(X,H)$ is similar to that of $k\sur B(X,H)$ and we have $k\sur A(X,H)\neq0$ if and only if $H$ is isomorphic to a subgroup of~$X$. Then the construction of the two kinds of bilinear forms on~$k\sur A(X,H)$ follow in the same manner and all the results of Sections~6 -- 8 hold. Note that the procedure for computing the dimension of $S_{H,V}(G)$ described at the end of Section~7 becomes easier, for it only involves the Mackey formula instead of its generalized version.\par

In this context of bifree bisets, the following new feature appears in characteristic zero.

\result{Proposition} \label{bifree}
Assume $k$ is a field of characteristic zero.
Let $X$ and $H$ be finite groups, and consider the symmetric bilinear form on~$k\sur A(X,H)$ defined by
$${<}\sur\alpha,\sur\beta{>}_X = \tau \pi(\alpha\op \beta) \,,$$
where $\alpha,\beta\in kA(X,H)$ are representatives of $\sur\alpha$ and $\sur\beta$ respectively, $\pi$ denotes the quotient map $\pi:kA(H,H)\to k\Out(H)$, and $\tau$ is the usual symmetrizing form on $k\Out(H)$ (coefficient of~1). Then this bilinear form is non-degenerate.
\fresult

\pf
The standard quotient $k\sur A(X,H)$ has a basis consisting of the elements $\Ind_A^X\, \Iso_\sigma$, where $A$ is a subgroup of~$X$ and $\sigma:H\to A$ is a group isomorphism. Moreover $A$ must be considered up to $X$-conjugation and $\Iso_\sigma$ up to right composition with inner automorphisms of~$H$ and left composition with conjugations by elements of $N_X(A)$. We compute the matrix of the bilinear form with respect to this basis and we show that it is a diagonal matrix.

Let $A$ and $B$ be subgroups of $X$ and let $\sigma:H\to A$ and $\rho:H\to B$ be group isomorphisms, considered up to conjugation as above. By the Mackey formula, we have
\begin{align*}
{<}\Ind_A^X\, \Iso_\sigma \,,\,& \Ind_B^X\, \Iso_\rho{>}_X
= \tau\pi(\Iso_\sigma^{-1}\,\Res_A^X\,\Ind_B^X\,\Iso_\rho) \\
&=\sum_{g\in[A\backslash X/B]}
\tau\pi(\Iso_\sigma^{-1}\,\Ind_{A\cap\ls gB}^A\,\Res_{A\cap\ls gB}^{\ls gB}\,\Conj_g\,\Iso_\rho) \,.
\end{align*}
This is zero if $B$ is not conjugate to $A$, because then $A\cap\ls gB$ is a proper subgroup of~$A$ and $\pi$ is zero on bisets which factorize through a proper subgroup of~$H$. Thus we assume now that $B$ is conjugate to $A$, and in fact $B=A$ without loss of generality. Then we obtain
\begin{align*}
{<}\Ind_A^X\, \Iso_\sigma \,,\,& \Ind_A^X\, \Iso_\rho{>}_X = \\
&=\sum_{g\in[A\backslash X/A]}
\,\tau\pi(\Iso_\sigma^{-1}\,\Ind_{A\cap\ls gA}^A\,\Res_{A\cap\ls gA}^{\ls gA}\,\Conj_g\,\Iso_\rho) \\
&=\sum_{g\in[N_X(A)/A]} \,\tau\pi(\Iso_\sigma^{-1}\,\Conj_g\,\Iso_\rho) \\
&=\sum_{g\in[N_X(A)/AC_X(A)]} \,|AC_X(A):A| \, \tau\pi(\Iso_\sigma^{-1}\,\Conj_g\,\Iso_\rho) \,.
\end{align*}
The term indexed by~$g$ is zero if $\pi(\sigma^{-1}\,\Conj_g\,\rho)$ is not the identity element of $\Out(H)$, so at most one term of the sum is non-zero. More precisely, the sum is zero unless there exists $g\in N_X(A)$ such that $\pi(\sigma^{-1}\, \Conj_g\,\rho)=1$ in $\Out(H)$, that is, $\sigma=\Conj_g\,\rho \,\Conj_h$ for some $h\in H$. But this means that we have the same basis element
$$\Ind_A^X\, \Iso_\sigma =\Ind_A^X\, \Conj_g\,\rho \,\Conj_h=\Ind_A^X\, \Iso_\rho \,.$$
This shows that the matrix is diagonal. Moreover, for every basis element $\Ind_A^X\, \Iso_\sigma$, the corresponding diagonal entry is $|AC_X(A):A|$. Since $k$ has characteristic zero, the bilinear form is non-degenerate.
\endpf

The non-degeneracy of the bilinear form implies the following special case of results of Webb~\cite{We2}.

\result{Corollary} \label{bifree-semi-simple}
Assume $k$ is a field of characteristic zero.
\begin{enumerate}
\item For every finite group $H$, the biset functor $k\sur A(-,H)$ is semi-simple.
\item For every finite groups $H$ and $G$, the $kA(G,G)$-module $k\sur A(G,H)$ is semi-simple.
\item The $k$-algebra $kA(G,G)$ is semi-simple.
\end{enumerate}
\fresult

\pf (1) Since the kernel of the bilinear form on~$k\sur A(X,H)$ is $0$ for each $X$, the functor $k\sur A(-,H)$ is semi-simple (Corollary~\ref{E-semi-simple}).\par

(2) This follows form (1) and evaluation at $G$.\par

(3) We only sketch the argument. It is not hard to prove that $kA(G,G)$ acts faithfully on $\prod_{H\leq G} k\sur A(G,H)$. Since this module is semi-simple by~(2), it follows that the radical $J(kA(G,G))$ must be zero.
\endpf

This corollary is a special case of more general results of Webb, who proved that the category of all bifree biset functors is semi-simple in characteristic zero (see Theorem~9.5 in~\cite{We2}). In particular the algebra $kA(G,G)$ is semi-simple, because it is the endomorphism algebra of the representable functor $kA(-,G)$. So we see that Corollary~\ref{bifree-semi-simple} is a special case of Webb's results.

Also, note that the question of vanishing evaluations, which is not easy in the general case of simple biset functors (see~\cite{BST}), has a direct solution for the category of global Mackey functors. Indeed, Webb proved an explicit formula giving the evaluation of simple functors (see Theorem~2.6 in~\cite{We1}).


\section{Examples}
\noindent
We illustrate the results of this paper by a few examples, which have been worked out either by hand or by computer calculations using~\cite{GAP}. We first start with a very small example, where most computations can be made over~$\Z$.

\result{Example} \label{Cp} 
{\rm 
Let $p$ be a prime and let $G=C_p$ be a cyclic group of order~$p$.
Then $\sur B(G,G)=\Z\Out(G)\cong \Z C_{p-1}$ and there is a ring homomorphism
$$\zeta: B(G,G) \longrightarrow \Z\Out(G) \cong \Z C_{p-1} \,.$$
On the other hand $\sur B(G,1)=B(G,1)$ is free abelian on the two elements $\Inf_{G/G}^G$ and $\Ind_1^G$. The action of $B(G,G)$ on $B(G,1)$ is easy to describe and this yields a ring homomorphism
$$\eta: B(G,G)\longrightarrow \End_{\Z}(B(G,1))\cong M_2(\Z)$$
whose image is generated by $\begin{pmatrix}1&0\\ 0&1 \end{pmatrix}$, $\begin{pmatrix}1&1\\ 0&0 \end{pmatrix}$, $\begin{pmatrix}1&p\\ 0&0 \end{pmatrix}$, $\begin{pmatrix}0&0\\ 1&1 \end{pmatrix}$, and $\begin{pmatrix}0&0\\ 1&p \end{pmatrix}$. Then we obtain an injective ring homomorphism
$$\zeta \times \eta: B(G,G) \longrightarrow \Z C_{p-1}\times M_2(\Z)$$
which provides an explicit representation of $B(G,G)$.\par

Extending scalars to a field $k$ whose characteristic does not divide $p-1$, we get an isomorphism
$$\zeta \times \eta: kB(G,G) \flh{\sim}{} kC_{p-1}\times M_2(k)$$
and we see that $kB(G,G)$ is semi-simple. If $k$ contains $(p-1)$-th roots of unity, then $kC_{p-1}$ decomposes further as a product of copies of~$k$, so $kB(G,G)$ has $p-1$ simple modules of dimension~1 and one simple module of dimension~2.
}
\fresult

In our next two examples, the double Burnside ring is not semi-simple, but we have a clear description of the Jacobson radical.

\result{Example} \label{D8} 
{\rm 
Let $D_8$ be the dihedral group of order~8 and let $k$ be a field of characteristic different from~2 and~3. By using Theorem~\ref{dimension-simple} (or its forerunner, Proposition~4.4.6 in~\cite{Bo3}, in the case of the trivial module~$k$), the dimension of $S_{H,V}(D_8)$ can be easily computed. Writing $(H,V)$ in the first line, we get the following values for $\dim(S_{H,V}(D_8))$:
$$\begin{array}{ccccccccc}
(1,k) & (C_2,k) & (C_4,k) & (C_4,\varepsilon) & (V_4,k) & (V_4, \varepsilon) & (V_4,2) & (D_8,k) & (D_8,\varepsilon) \\
5& 11 & 1 & 0 & 3 & 1 & 4 & 1 & 1
\end{array}$$
Here $\varepsilon$ denotes the sign representation of $\Out(H)$, the group $V_4$ denotes the Klein four group, and 2 denotes the two-dimensional representation of $\Out(V_4)$ (which is the symmetric group of order~6). Summing up the squares of the dimensions, we find $\dim(kB(D_8,D_8)/J)=175$, where $J=J(kB(D_8,D_8))$. On the other hand, the basis of Lemma~\ref{biset} yields $\dim(kB(D_8,D_8))=214$, hence $\dim(J)=214-175=39$. But we know, by Theorem~\ref{trivial-group}, that the ideal $I(G)\cap J$ has dimension $b(D_8)^2-c(D_8)^2=8^2-5^2=39$. It follows that $I(G)\cap J=J$, so the only contribution to the Jacobson radical comes from the trivial group.
}
\fresult

\result{Example} \label{Q8} 
{\rm 
Let $Q_8$ be the quaternion group of order~8 and let $k$ be a field of characteristic different from~2 and~3. This example is similar to the previous one and the computations show again that $I(G)\cap J=J$, so the only contribution to the Jacobson radical comes from the trivial group.
}
\fresult

We now move to a slightly more involved case.

\result{Example} \label{A4-C3}
{\rm 
Let $G=A_4$ and $H=C_3$ and suppose that $k$ is a field of characteristic different from~2. Then $k\sur B(A_4,C_3)$ has dimension~4, its radical $J$ has dimension~2 and coincides with the kernel $R(A_4,C_3)$ of the bilinear form (by Proposition~\ref{abelian-etc}). Moreover
$$k\sur B(A_4,C_3)/J=S_{C_3,k_+}(A_4)\oplus S_{C_3,k_-}(A_4)
\quad\text{and}\quad J=S_{A_4,k_+}(A_4)\oplus S_{A_4,k_-}(A_4) \,,$$
where $k_+$ denotes the trivial representation of the cyclic group $\Out(C_3)$ of order~2 and $k_-$ denotes the sign representation. All the evaluations above are 1-dimensional. This also shows, 
not surprisingly, that the Jacobson radical of $kB(G,G)$ is not entirely contained in the ideal $I(G)$ considered in Section~10.
}
\fresult

Our next examples are concerned with the fact that $J\big(k\sur B(G,H)\big)$ is not necessarily equal to $R(G,H)$, though many cases where equality holds have been seen in Section~9.

\result{Example} \label{example-A5}
{\rm 
Let $G=A_5$ be the alternating group on 5 letters, let $H=C_3$ be the cyclic group of order~3, and suppose that $k$ is a field of characteristic different from~2. The vector space $k\sur B(A_5,C_3)$ has dimension~3, with basis $\Ind_{C_3}^{A_5}, \Indinf_{A_4/V_4}^{A_5}, \Indinf_{A_4/V_4}^{A_5} \Iso_\sigma$, where $V_4$ is the Klein four-group and $\sigma$ is the non-trivial group automorphism of the cyclic group of order~3 (see Lemma~\ref{basis} and note that $\Ind_{C_3}^{A_5} \Iso_\sigma \cong \Ind_{C_3}^{A_5}$.)\par

We have $E=k\Out(C_3)\cong k_+\times k_-$, so there are two simple functors $S_{C_3,k_+}$ and $S_{C_3,k_-}$ indexed by~$C_3$. The direct computation of the bilinear form shows that $R(A_5,C_3)$ has dimension~2, generated by the differences of basis elements, so $k\sur B(A_5,C_3)/R(A_5,C_3)$ is one-dimensional. We have in fact
$$k\sur B(A_5,C_3)/R(A_5,C_3) \cong S_{C_3,k_+}(A_5) \oplus S_{C_3,k_-}(A_5) = S_{C_3,k_+}(A_5) = k \,,$$
because actually $S_{C_3,k_+}(A_5) = k$ and $S_{C_3,k_-}(A_5) = 0$ in this specific example.\par

But further computations show that $J\big(k\sur B(A_5,C_3)\big)$ is only one-dimensional and that we have to take into account the simple functors $S_{A_4,V}$ indexed by the larger group $A_4$. We have to consider the simple $kB(A_5,A_5)$-modules $S_{A_4,V}(A_5)$, where $V=k_+$ or $V=k_-$, the two simple modules for the cyclic group $\Out(A_4)$ of order~2. It turns out that we get a semi-simple quotient of dimension~2
$$k\sur B(A_5,C_3)/J\big(k\sur B(A_5,C_3)\big) \cong S_{C_3,k_+}(A_5) \oplus S_{A_4,k_-}(A_5) \,,$$
with an extra factor indexed by~$A_4$. Moreover, the Jacobson radical has dimension~1 and satisfies $J\big(k\sur B(A_5,C_3)\big)\cong S_{A_4,k_+}(A_5)$.
}
\fresult

\result{Example} \label{other-examples}
{\rm 
There are numerous other examples where the radical $J=J\big(k\sur B(G,H)\big)$ of the module $M=k\sur B(G,H)$ is not equal to the kernel $R=R(G,H)$ of the bilinear form. In such cases, there are additional simple quotients, namely the factors of $R/J$, which are indexed by groups larger than~$H$. We just list a few such examples.\par

$G=GL(3,2)$ and $H=C_3$. Then $\dim(M)=5$, $\dim(R)=4$, $\dim(J)=3$.\par

$G=SL(2,7)$ and $H=C_3$. Then $\dim(M)=8$, $\dim(R)=6$, $\dim(J)=4$. Here, $M/R=S_{C_3,k}(G)$, with dimension~2, but $S_{C_3,k_-}(G)=0$. There is one additional simple quotient $R/J$, of dimension~2, indexed by the group $C_7\rtimes C_3$. Moreover, $J$ is the direct sum of two 2-dimensional simple modules, one indexed by $A_4$, and the other by $C_7\rtimes C_3$.\par

$G=\tilde{A}_5$ (perfect group of order 120) and $H=C_3$. Then $\dim(M)=4$, $\dim(R)=2$, $\dim(J)=1$. The additional simple quotient $R/J$ turns out to be indexed by the group $A_4$.\par

$G=\tilde{A}_5$ and $H=C_4$. Then $\dim(M)=5$, $\dim(R)=4$, $\dim(J)=2$. Here, $M/R=S_{C_4,k}(G)$, with dimension~1, but $S_{C_4,k_-}(G)=0$. There are 2 additional simple quotients appearing in $R/J$, one indexed by the group $C_5\rtimes C_4$, and the other by the group $C_3\rtimes C_4$.\par

$G=PSL(2,11)$ and $H=C_3$. Then $\dim(M)=4$, $\dim(R)=2$, $\dim(J)=1$. The additional simple quotient $R/J$ turns out to be indexed by the group $A_4$.\par

$G=PSL(2,11)$ and $H=C_5$. Then $\dim(M)=6$, $\dim(R)=4$, $\dim(J)=2$.
Here, the evaluations $S_{C_5,V}(G)$ vanish for two of the $k\Out(C_5)$-modules~$V$, while the other two appear in the quotient~$M/R$.\par

$G=PSL(2,8)$ and $H=C_7$. Then $\dim(M)=9$, $\dim(R)=6$, $\dim(J)=3$. There are 3 additional simple quotients appearing in $R/J$, indexed by the group $(C_2)^3\rtimes C_7$ (the normalizer of a Sylow 2-subgroup of $G$).\par
}
\fresult

\bigskip
\noindent
Serge Bouc, CNRS-LAMFA, Universit\'e de Picardie - Jules Verne,\\
33, rue St Leu, F-80039 Amiens Cedex~1, France.\\
{\tt serge.bouc@u-picardie.fr}

\medskip
\noindent
Radu Stancu, CNRS-LAMFA, Universit\'e de Picardie - Jules Verne,\\
33, rue St Leu, F-80039 Amiens Cedex~1, France.\\
{\tt radu.stancu@u-picardie.fr}

\medskip
\noindent
Jacques Th\'evenaz, Section de math\'ematiques, EPFL, \\
Station~8, CH-1015 Lausanne, Switzerland.\\
{\tt Jacques.Thevenaz@epfl.ch}

\end{document}